\documentclass{gtart_a}
\pdfoutput=1
\usepackage{pinlabel}

\title{The canonical pencils on Horikawa surfaces}

\author{Denis Auroux}
\givenname{Denis}
\surname{Auroux}
\address{Department of Mathematics\\MIT\\Cambridge, MA 02139\\USA}
\email{auroux@math.mit.edu}
\urladdr{}

\volumenumber{10}
\issuenumber{}
\publicationyear{2006}
\papernumber{47}
\lognumber{0747}
\startpage{2173}
\endpage{2217}

\doi{}
\MR{}
\Zbl{}


\received{28 May 2006}
\revised{25 September 2006}
\accepted{22 October 2006}
\proposed{Rob Kirby}
\seconded{Ronald Fintushel, Ronald Stern}
\published{29 November 2006}
\publishedonline{29 November 2006}
\corresponding{}
\editor{CPR}
\version{}

\arxivreference{math.GT/0605692}

\let\xysavmatrix\xymatrix
\def\xymatrix{\disablesubscriptcorrection\xysavmatrix}
\AtBeginDocument{\let\bar\wbar\let\tilde\wtilde\let\hat\what}

\makeatletter
\def\cnewtheorem#1[#2]#3{\newtheorem{#1}{#3}[section]
\expandafter\let\csname c@#1\endcsname\c@theorem}

\theoremstyle{plain}
\newtheorem{theorem}{Theorem}[section]
\cnewtheorem{lemma}[theorem]{Lemma}
\cnewtheorem{proposition}[theorem]{Proposition}
\cnewtheorem{corollary}[theorem]{Corollary}
\cnewtheorem{conj}[theorem]{Conjecture}
\cnewtheorem{question}[theorem]{Question}
\theoremstyle{definition}
\cnewtheorem{definition}[theorem]{Definition}
\cnewtheorem{remark}[theorem]{Remark}

\makeatother  

\newcommand{\CP}{\mathbb{CP}}
\newcommand{\FF}{\mathbb{F}}
\newcommand{\PP}{\mathbb{P}}
\newcommand{\OO}{\mathcal{O}}
\newcommand{\DD}{\mathsf{D}}
\newcommand{\crit}{\mathrm{crit}}
\newcommand{\Map}{\mathrm{Map}{}}
\newcommand{\unsp}{\hspace{-0.06em}}
\newcommand{\JJ}{\mathbb{J}}

\begin{document}

\begin{abstract}
We calculate the monodromies of the canonical Lefschetz pencils on a
pair of homeomorphic Horikawa surfaces. We show in particular that the 
(pluri)canonical pencils on these surfaces have the same monodromy
groups, and are related by a ``partial twisting'' operation.
\end{abstract}

\maketitle

\section{Introduction}

Horikawa surfaces are minimal complex surfaces of general type which
realize the equality case in Noether's inequality $c_1^2\ge 2p_g-4$.
While their classification as complex surfaces has been completed a
long time ago \cite{Hor},
the topology of these surfaces viewed as smooth 4--manifolds, or as symplectic
4--manifolds, remains mysterious. In this paper we consider two specific
Horikawa surfaces:

\begin{definition}
Denote by $X_1$ a double cover of\/ $\CP^1\times\CP^1$ branched along
a smooth algebraic curve $C_1$ of bidegree $(6,12)$. Denote by $X_2$ a
double cover of the Hirzebruch surface\/ $\FF_6=\PP(\OO_{\PP^1}\oplus
\OO_{\PP^1}(6))$ branched along $\Delta_\infty\cup C_2$, where $\Delta_\infty$
is the exceptional section of\/ $\FF_6$ $(\Delta_\infty\cdot\Delta_\infty=-6)$,
and $C_2$ is a smooth algebraic curve in the linear system $|5\Delta_0|$
(where $\Delta_0$ is the zero section, satisfying $\Delta_0\cdot \Delta_0=+6$).
\end{definition}

\noindent
(The actual choices of $C_1$ and $C_2$ are irrelevant from the point of view
of symplectic topology, hence we do not specify them here).

The complex surfaces $X_1$ and $X_2$ are simply-connected, non-spin, and
have the same Euler characteristic $e(X_1)=e(X_2)=116$ and signature
$\sigma(X_1)=\sigma(X_2)=-72$; hence by a classical result of Freedman
they are homeomorphic. Moreover, the homeomorphism between them can be
chosen so that the canonical classes $K_{X_1}$ and $K_{X_2}$ are the image
of each other under the induced map on cohomology. However, $X_1$ and $X_2$
are not deformation equivalent as complex surfaces \cite{Hor}. The question
of whether $X_1$ and $X_2$ are diffeomorphic is open to this date; although
the expected answer is negative, a result of Friedman and Morgan \cite{FM}
shows that $X_1$ and $X_2$ cannot be distinguished using Seiberg--Witten theory.

Since the canonical classes $K_{X_1}$, $K_{X_2}$ are ample, one can equip
$X_1$ and $X_2$ with K\"ahler forms such that $[\omega_i]=c_1(K_{X_i})$ and
view them as symplectic 4--manifolds. The question of whether $(X_1,\omega_1)$
and $(X_2,\omega_2)$ are symplectomorphic is again open; following a
strategy proposed by Donaldson \cite{Do1}, one can try to approach this
problem from the perspective of Lefschetz pencils.

Namely, considering a generic pencil of curves in the linear system 
$|K_{X_i}|$ and blowing up its 16 base points, we obtain a fibration
$\hat{f}_i\co \hat{X}_i\to\CP^1$. The generic fiber of $\hat{f}_i$ is a smooth curve
of genus 17, and the singular fibers are nodal (there are 196 nodes in
total); moreover, the exceptional
divisors of the blowups determine 16 distinguished sections of $\hat{f}_i$.
It is well-known (see \fullref{ss:slf}) that the monodromy of such a 
fibration can be described by an ordered tuple of Dehn twists in the
mapping class group $\Map_{17,16}$ of a genus 17 surface with 16 boundary
components (up to Hurwitz equivalence and global conjugation).

\begin{theorem} \label{thm:main1}
The canonical pencils on $X_1$ and $X_2$ are related by a ``partial twisting''
operation, ie, there exist Dehn twists $\phi,t_1,\dots,t_{196}
\in\Map_{17,16}$ such that the monodromy of $\hat{f}_1$ can be expressed by
the tuple $(\phi t_1\phi^{-1},\dots,$ $\phi t_{64}\phi^{-1}, t_{65},
\dots, t_{196})$, and the monodromy of $\hat{f}_2$ can be expressed by the
tuple $(t_1,\dots,t_{64},t_{65},\dots,t_{196})$.
\end{theorem}

There is a geometric reason for this property of the canonical pencils:
the Horikawa surfaces $X_1$ and $X_2$ can be obtained
from each other by a {\it Luttinger surgery} operation
(cf \fullref{thm:horilu}).

In fact, we determine the monodromies of $\hat{f}_1$ and $\hat{f}_2$ explicitly
(the formulas are given in Theorems \ref{thm:monx2} and \ref{thm:monx1});
as a consequence, we also get:

\begin{theorem} \label{thm:subgp}
The monodromy groups of $\hat{f}_1$ and $\hat{f}_2$, ie, the subgroups of
$\Map_{17,16}$ generated by the Dehn twists in their monodromy, are
isomorphic to each other.
\end{theorem}

All these results suggest that (somewhat unsurprisingly) the pencils $\hat{f}_1$ and
$\hat{f}_2$ are very similar to each other, and difficult to tell apart. This is in sharp
contrast with the well-known genus 2 fibrations carried by $X_1$ and $X_2$,
whose monodromies are easily distinguished (see \fullref{ss:horig2}).

By the work of Gompf \cite{GS,Go}, if the tuples of Dehn twists
describing $\hat{f}_1$ and $\hat{f}_2$ are Hurwitz and conjugation equivalent, then
$(X_1,\omega_1)$ and $(X_2,\omega_2)$ are symplectomorphic. The converse is
not necessarily true. However, by Donaldson's asymptotic uniqueness result
for symplectic
Lefschetz pencils \cite{Do2}, if $(X_1,\omega_1)$ and
$(X_2,\omega_2)$ are symplectomorphic then there exists an
integer $k_0$ such that the pluricanonical Lefschetz pencils on $X_1$ and
$X_2$ (ie, generic pencils of curves in the linear systems
$|kK_{X_i}|$) have equivalent monodromies for all $k\ge k_0$.
Hence, in order to prove that $(X_1,\omega_1)$ and $(X_2,\omega_2)$ are
not symplectomorphic, one needs to compare not just the canonical pencils
of $X_1$ and $X_2$, but also a sequence of pluricanonical pencils.

``Degree doubling'' arguments \cite{AK,Sm} (see also \fullref{ss:doubling}) imply that the
monodromies of the pluricanonical pencils (for the linear systems
$|2^mK_{X_i}|$) are determined by those of the canonical pencils in an
explicit and
``universal'' manner. In particular, it is expected that an invariant that
distinguishes the canonical pencils should also be able to distinguish
the pluricanonical pencils (and hence prove that $X_1$ and $X_2$ are not
symplectomorphic). While no such invariant is known to this date,
it is still worth mentioning the following consequence of
Theorems \ref{thm:main1} and \ref{thm:subgp}:

\begin{theorem}\label{thm:pluri}
For all $m\ge 0$, generic pencils of curves in the linear systems $|2^m
K_{X_i}|$ on $X_1$ and $X_2$ are related to each other by ``partial twisting'' operations,
and their monodromy subgroups are isomorphic.
\end{theorem}

The rest of this paper is organized as follows: in \fullref{s:background}
we review some background material on Lefschetz fibrations and their monodromy
(\fullref{ss:slf}), lifting homomorphisms (\fullref{ss:lifting}), Horikawa
surfaces (\fullref{ss:horig2}), and Luttinger surgery (\fullref{ss:luttinger}).
\fullref{s:monx2} is devoted to the calculation of the monodromy of the canonical pencil
on $X_2$. In \fullref{s:monx1} we show that $X_1$ and $X_2$ are related
by Luttinger surgery, and which allows us to determine the
monodromy of a certain symplectic Lefschetz pencil on $X_1$;
using the theory of pseudo-holomorphic curves, we prove in \fullref{s:jhol}
that this Lefschetz pencil is isomorphic to the canonical pencil
of $X_1$. This allows us to complete the proof of \fullref{thm:main1},
while \fullref{thm:subgp} is proved in \fullref{s:mainthm}. 
The paper ends with considerations about
degree doubling and pluricanonical pencils in \fullref{ss:doubling},
and Lagrangian spheres and matching paths in \fullref{ss:matching}.

While we are still a long way from proving that $X_1$ and $X_2$ are not
symplectomorphic, the explicit calculation of the monodromies of their
canonical pencils sheds some light on the situation; we hope that it
will lead to further advances on this problem, and give some
insight about what kind of invariants one might consider in order to 
distinguish homeomorphic surfaces of general type.
\medskip

 \noindent
{\bf Acknowledgements}\qua This work was partially supported by NSF
grants DMS-0244844 and DMS-0600148 and an A\,P Sloan research
fellowship. The author wishes to thank the Department of Mathematics
at UC Berkeley for its hospitality during part of the preparation of
this work.

\section{Preliminaries} \label{s:background}

\subsection{Lefschetz fibrations and symplectic 4--manifolds}\label{ss:slf}

\begin{definition}
A Lefschetz fibration on an oriented compact smooth 4--mani\-fold $M$ is a
smooth map $f\co M\to S^2$ which is a submersion everywhere except at finitely
many non-degenerate critical points $p_1,\dots,p_r$, near which $f$
identifies in local orientation-preserving complex coordinates with
the model map $(z_1,z_2)\mapsto z_1^2+z_2^2$.
\end{definition}

The fibers of a Lefschetz fibration
$f$ are compact oriented surfaces, smooth except for finitely
many of them. The fiber through $p_i$ presents a transverse double point,
or {\it node}, at $p_i$. Without loss of generality, we can assume after
perturbing $f$ slightly that the critical values $q_i=f(p_i)$ are all distinct.
Fix a reference point $q_*$ in $S^2\setminus \crit(f)$, and let
$\Sigma=f^{-1}(q_*)$ be the corresponding fiber. Then we can consider
the {\it monodromy homomorphism}
$$\psi\co \pi_1(S^2\setminus\crit(f),q_*)\to \Map(\Sigma),$$
where $\Map(\Sigma)=\pi_0\mathrm{Diff}^+(\Sigma)$ is the mapping class group
of $\Sigma$. The image $\psi(\gamma)$ of a loop
$\gamma\subset S^2\setminus\crit(f)$ is
the isotopy class of the diffeomorphism of $\Sigma$ induced by parallel
transport (with respect to an arbitrary horizontal distribution)
along the loop $\gamma$.

The singular fibers of $f$ are obtained from the nearby smooth fibers by
collapsing a simple closed loop, called the {\it vanishing cycle}.
The monodromy of a Lefschetz fibration around a singular fiber is
the positive Dehn twist along the corresponding vanishing cycle.
Choose an ordered collection $\eta_1,\dots,\eta_r$ of arcs joining $q_*$
to the various critical values of $f$, and thicken them to obtain closed
loops $\gamma_1,\dots,\gamma_r$ based at $q_*$ in
$S^2\setminus \crit(f)$, such that each $\gamma_i$ encircles exactly one
of the critical values of $f$, and
$\pi_1(S^2\setminus \crit(f),q_*)=\langle \gamma_1,\dots,\gamma_r\,|\,
\prod \gamma_i=1\rangle.$ Then the monodromy of $f$ along each $\gamma_i$ is
a positive Dehn twist $t_i$ along an embedded loop
$\delta_i\subset\Sigma$, obtained by parallel transport along $\eta_i$
of the vanishing
cycle at the critical point $p_i$, and in $\Map(\Sigma)$ we have the relation
$t_1\dots t_r=\mathrm{Id}.$

Hence, to every Lefschetz fibration we can associate a {\it factorization} of
the identity element as a product of positive Dehn twists in the mapping class
group of the fiber, ie, an ordered tuple of
Dehn twists whose product is equal to Id; we will often use the multiplicative
notation, with the understanding that what is important
is not the product of the factors but rather the factors themselves.

Given the collection of Dehn twists $t_1,\dots,t_r$
we can reconstruct the Lefschetz fibration $f$ above a large disc $\DD$
containing all the critical values, by starting from $\Sigma\times D^2$ and
adding handles as specified by the vanishing cycles \cite{Kas}. To recover
the 4--manifold $M$ we need to glue $f^{-1}(\DD)$ and the trivial
fibration $f^{-1}(S^2\setminus\DD)=\Sigma\times D^2$ along their
common boundary, in a
manner compatible with the fibration structures. In general this gluing
involves the choice of an element in $\pi_1\mathrm{Diff}^+(\Sigma)$; however
the diffeomorphism group is simply connected if the genus of $\Sigma$ is at
least 2, and in that case the factorization $t_1\dots t_r=\mathrm{Id}$
determines the Lefschetz fibration $f\co M\to S^2$ completely (up to isotopy).

The monodromy factorization $t_1\dots t_r=\mathrm{Id}$
depends not only on the topology of $f$,
but also on the choice of an ordered collection $\gamma_1,\dots,\gamma_r$
of generators of $\pi_1(S^2\setminus \crit(f),q_*)$; the braid group
$B_r$ acts transitively on the set of all such ordered collections, by
{\it Hurwitz moves}. The equivalence relation induced by this action on
the set of mapping class group factorizations is generated by
$$(t_1,\dots,t_i,t_{i+1},\dots,t_r)\sim (t_1,\dots,
t_it_{i+1}t_i^{-1},t_i,\dots,t_r)\quad\forall 1\le i<r,$$
and is called {\it Hurwitz equivalence}. Additionally, in order to remove
the dependence on the choice of the reference fiber $\Sigma$, we should
view the Dehn twists $t_i$ as elements of the mapping class group
$\Map_g$ of an abstract surface of genus $g=g(\Sigma)$. This requires the
choice of an identification diffeomorphism, and introduces another
equivalence relation on the set of mapping class group factorizations:
{\it global conjugation},
$$(t_1,\dots,t_r)\sim (\phi t_1 \phi^{-1},\dots,\phi t_r
\phi^{-1})\quad \forall \phi\in\Map_g.$$

\begin{proposition}\label{prop:mono}
For $g\ge 2$, there is a one to one correspondence between $(a)$
factorizations of \,{\rm Id} as a product of positive Dehn twists in
$\Map_g$, up to Hurwitz equivalence and global conjugation, and $(b)$
genus $g$ Lefschetz fibrations over $S^2$, up to isomorphism.
\end{proposition}

It is a classical result of Thurston that, if $M$ is an oriented
surface bundle over an oriented surface, then $M$ is a symplectic
4--manifold, at least provided that the homology class of the fiber
is nonzero in $H_2(M,\R)$. As shown by Gompf, the argument extends
to the case of Lefschetz fibrations  \cite[Theorem 10.2.18]{GS}:

\begin{theorem}[Gompf]\label{thm:go}
Let $f\co M\to S^2$ be a Lefschetz fibration, and assume that the fiber
represents a nonzero class in $H_2(M,\R)$. Then $M$ admits a symplectic
structure for which the fibers of $f$ are symplectic submanifolds; this
symplectic structure is unique up to deformation.
\end{theorem}

Lefschetz fibrations arise naturally in algebraic geometry: 
if $X$ is a complex surface, and $L\to X$ is a sufficiently ample line
bundle, then the ratio between two suitably chosen sections $s_0,s_1\in
H^0(L)$ determines a {\it Lefschetz pencil}, ie, a map
$f=(s_0/s_1)\co X\setminus B\to\CP^1$, defined on the complement of the finite
set $B=\{s_0=s_1=0\}$ (the {\it base points}), with isolated nondegenerate
critical points in $X\setminus B$. 

More generally, Donaldson has shown that this construction extends to the
symplectic setting \cite{Do2}:

\begin{theorem}[Donaldson]\label{thm:do}
Let $(X,\omega)$ be a compact symplectic 4--manifold with
$[\omega]\in H^2(X,\Z)$. Then $X$ carries a {\em
symplectic Lefschetz pencil}, ie, there exist
a finite set $B\subset X$ and a map $f\co X\setminus B\to
\CP^1=S^2$ such that $f$ is modelled on $(z_1,z_2)\mapsto (z_1:z_2)$ near
each point of $B$, and $f$ is a Lefschetz fibration with
symplectic fibers outside of $B$.
\end{theorem}

More precisely, given a compact symplectic 4--manifold
$(X,\omega)$ such that $[\omega]\in H^2(X,\Z)$, and given a complex
line bundle $L\to X$ with $c_1(L)=[\omega]$, one can construct symplectic
Lefschetz pencils from suitably chosen pairs of sections $s_0,s_1\in
C^\infty(L^{\otimes k})$ for all sufficiently large values of the integer $k$.
Moreover, Donaldson has shown that for all sufficiently large values 
$k$, there is a distinguished connected component of the space of
symplectic Lefschetz pencils obtained from pairs of sections of
$L^{\otimes k}$; in the K\"ahler case, this component contains pencils defined
from pairs of generic holomorphic sections of $L^{\otimes k}$ \cite{Do2}.
Hence, if two complex projective surfaces are symplectomorphic, then generic
pencils of curves in the linear systems considered by Donaldson are mutually
isomorphic (as symplectic Lefschetz pencils) whenever the integer $k$ 
is sufficiently large.

Given a symplectic Lefschez pencil $f\co X\setminus B\to S^2$, the
manifold $\hat{X}$ obtained from $X$ by
blowing up the points of $B$ admits a Lefschetz fibration
$\hat{f}\co \hat{X}\to S^2$ with symplectic fibers, and can be described by its
monodromy as discussed above.

Moreover, the fibration $\hat{f}$ has $n=|B|$ distinguished sections
$e_1,\dots,e_n$, corresponding to the exceptional divisors of the blowups.
Therefore, each fiber of $\hat{f}$ comes equipped with $n$ marked points,
and the monodromy of $\hat{f}$ lifts to the mapping class group of a genus
$g$ surface with $n$ marked points.

The normal bundles of the sections $e_i$ are not trivial, but
it is possible to trivialize them over the preimage of a
large disc $\DD$ containing all the chosen generators of
$\pi_1(S^2\setminus\crit(\hat{f}))$. Deleting a small tubular
neighborhood of each exceptional section, we can now view the monodromy of
$\hat{f}$ as a morphism
$$\hat\psi\co \pi_1(\DD\setminus \crit(\hat{f}))\to \Map_{g,n},$$
where $\Map_{g,n}$ is the mapping class group of a genus $g$ surface with
$n$ boundary components (ie, $\pi_0\mathrm{Diff}^+(\Sigma,\partial\Sigma)$).

The product of the Dehn twists
$t_i=\hat\psi(\gamma_i)$ is no longer the identity element in
$\Map_{g,n}$. Instead, since $\prod \gamma_i$ is homotopic to the boundary
of the disc $\DD$, and since the normal bundle to $e_i$ has degree $-1$,
we have $\prod t_i=T_\partial$, where $T_\partial\in\Map_{g,n}$ is the
{\it boundary twist}, ie, the product of the positive Dehn twists
along $n$ loops parallel to the boundary components.

With this understood, the previous discussion carries over, and under the
assumption $2-2g-n<0$ there is a
one to one correspondence between factorizations of the boundary twist
$T_\partial$ as a product of positive Dehn twists in $\Map_{g,n}$, up to Hurwitz
equivalence and global conjugation, and genus $g$ Lefschetz fibrations over
$S^2$ equipped with $n$ distinguished sections of square $-1$, up to
isomorphism.

Moreover, \fullref{thm:go} admits a strengthening in this context:
for symplectic Lefschetz pencils whose fibers are Poincar\'e dual to a
symplectic form, the monodromy data determines the
symplectic structure up to isotopy (ie, symplectomorphism), rather than
just up to deformation \cite{Go}. Combining this with the discussion after
\fullref{thm:do}, we conclude:

\begin{corollary}[Donaldson, Gompf]
The following three properties are equivalent:

\begin{itemize}
\item[{\rm(i)}] the Horikawa surfaces $X_1$ and $X_2$ equipped with their canonical
K\"ahler forms are symplectomorphic;

\item[{\rm(ii)}] there exists an integer $k\ge 1$ such that generic pencils 
of curves in the linear systems $|kK_{X_i}|$ have equivalent monodromy
factorizations;

\item[{\rm(iii)}] there exists an integer $k_0$ such that, for all $k\ge k_0$,
generic pencils of curves in the linear systems $|kK_{X_i}|$ have
equivalent monodromy factorizations.
\end{itemize}
\end{corollary}

\subsection{Double covers and lifting homomorphisms}\label{ss:lifting}

Let $X$ and $Y$ be smooth complex surfaces, such that there exists a 2:1
covering map $\pi\co X\to Y$, branched along a smooth curve $C\subset Y$.
Assume that we have a Lefschetz pencil $f\co Y\setminus B\to \CP^1$ (we also
allow $f$ to be a fibration, ie, $B$ may be empty), and that the branch
curve $C$ satisfies the following properties:%

\begin{itemize}
\item[{\rm(i)}] $C$ does not pass through the base
points or the critical points of $f$;

\item[{\rm(ii)}] $C$ is everywhere transverse to the fibers of $f$, except at isolated
points $p_1,\dots,p_s$ where $C$ is nondegenerately tangent to the fiber
of $f$ (ie, at $p_i$ the multiplicity of the intersection between $C$ and
the fiber of $f$ is 2);

\item[{\rm(iii)}] for simplicity we also assume that the
points $p_i$ lie in distinct smooth fibers of the pencil $f$.
\end{itemize}

Then $\tilde{f}=f\circ \pi\co X\setminus \tilde{B}\to\CP^1$ is also a Lefschetz
pencil, with base points $\tilde{B}=\pi^{-1}(B)$.

\begin{remark}
The discussion extends without modification to the situation
where $\pi\co X\to Y$ is a branched covering of symplectic $4$--manifolds with a
smooth symplectic branch curve, and $f$ is a symplectic Lefschetz pencil.
\end{remark}

Denote by $\Sigma$ the generic fiber
of $f$, with a neighborhood of the base points removed (so $\Sigma$ is a
compact surface with $n=|B|$ boundary components, and the monodromy of $f$
takes values in the mapping class group of $\Sigma$). The generic fiber of
$\tilde{f}$ is a double cover of $\Sigma$ branched at $d=[C]\cdot[\Sigma]$
points, which we denote by $\tilde{\Sigma}$. Abusing notation, we denote
the restriction of the double cover to the fiber by the same letter:
$\pi\co \tilde{\Sigma}\to\Sigma$. 

It is a classical fact that the double covering $\pi$ determines a lifting
homomorphism $L$ from (a subgroup of) the braid group $B_d(\Sigma)$
(ie, the fundamental group of the space
$\mathcal{C}_d(\Sigma)$ of unordered
configurations of $d$ distinct points in the interior of $\Sigma$)
to the mapping class group of $\tilde\Sigma$. In the case where $\Sigma$
has genus 0, which is the only one we will be considering, one
way to describe the lifting homomorphism is to consider the universal family
$\mathcal{X}\to \mathcal{C}_d(\Sigma)$ whose fiber above a configuration
$\{x_1,\dots,x_d\}\subset\Sigma$ is the double cover of $\Sigma$ branched
at $x_1,\dots,x_d$ (with trivial monodromy along each component of
$\partial\Sigma$). (When the genus of $\Sigma$ is nonzero, this universal
family is only defined over a finite covering of $\mathcal{C}_d(\Sigma)$).
The lifting homomorphism is simply the monodromy
of the fibration $\mathcal{X}\to\mathcal{C}_d(\Sigma)$; however not every element of $B_d(\Sigma)$ is liftable,
because some braids lift to diffeomorphisms of $\tilde\Sigma$ which exchange
the two lifts of some boundary components of $\Sigma$ instead of fixing $\partial
\tilde\Sigma$ pointwise.

Given any arc $\eta$ joining two points $x_i^0$ and $x_j^0$ of the
reference configuration\break
$\{x_1^0,\dots,x_d^0\}=\Sigma\cap C\subset\Sigma$ inside
$\Sigma\setminus \{x_1^0,\dots,x_d^0\}$, we can consider the {\it
half-twist} along $\eta$, which is the braid exchanging the two points 
$x_i^0$ and $x_j^0$ by a counterclockwise 180 degree rotation inside a
small tubular neighborhood of $\eta$. The preimage $\pi^{-1}(\eta)$ is a
simple closed curve in $\tilde\Sigma$, and it is a classical observation
that the half-twist along $\eta$ lifts to the Dehn twist along
$\pi^{-1}(\eta)$.%

We claim that the monodromy of the Lefschetz fibration $\tilde{f}$ is
completely determined by the monodromy of $f$ and by the {\it braid monodromy} of the
branch curve $C$. Indeed, the singular fibers of $\tilde{f}$ are of two types:

\begin{itemize}
\item[(a)] preimages by $\pi$ of the singular fibers of $f$;

\item[(b)] preimages by $\pi$ of smooth fibers of $f$ which are tangent to $C$.
\end{itemize}

Denote by $q_1,\dots,q_r$ the critical points of $f$, which we assume to
lie in distinct fibers, and choose an ordered collection of generating loops
for $\pi_1(\DD\setminus \crit(\tilde{f}))$, where $\DD$ is a large disc
containing all the points of $\crit(\tilde{f})=\{f(q_1),
\dots,f(q_r),$ $f(p_1),\dots,f(p_s)\}$. This allows us to define the monodromy
of $\tilde{f}$ around its various critical values as in \fullref{ss:slf}.

First consider a singular fiber of $f$, containing
a critical point $q_j$. The corresponding fiber of $\tilde{f}$ possesses
two nodal singularities (the two preimages of $q_j$), and we claim that the
corresponding vanishing cycles are the two lifts by $\pi$ of the vanishing
cycle at $q_j$.

Indeed, by assumption $q_j\not\in C$,
and since $C$ is transverse to the fibers of $f$ in the considered region,
we can locally choose the parallel transport maps between the
various fibers of $f$ in a manner such that the intersection points with $C$
are preserved. Therefore, the vanishing cycle of $f$ associated to the
critical point $q_j$ and to the chosen generator of $\pi_1(\DD\setminus
\crit(\tilde{f}))$ can be naturally represented by a simple closed curve
$\delta_j\subset\Sigma\setminus \{x_1^0,\dots,x_d^0\}$. 
The preimage $\pi^{-1}(\delta_j)$ consists of two disjoint simple closed
curves \smash{$\delta'_j$ and $\delta''_j$}, which are precisely the two
vanishing cycles of $\tilde{f}$. The monodromy of $f$ along the chosen
loop around $q_j$ is the Dehn
twist along $\delta_j$, and the monodromy of $\tilde{f}$ along the same loop
is the product of the Dehn twists along $\delta'_j$ and~$\delta''_j$.

\begin{remark}
By construction, the critical values of $\tilde{f}$ are not distinct,
since critical points of $f$ lift to pairs of critical points in the
same fiber. However, when the chosen linear system is sufficiently ample,
the critical values can be made distinct by considering a small generic
perturbation of $\tilde{f}$; if the perturbation is chosen
sufficiently small then the vanishing cycles are not affected. More
generally, even when such a perturbation does not exist in the algebraic
setting (as is the case for canonical pencils on Horikawa surfaces), it
can still be carried out among symplectic Lefschetz pencils. Hence, when
viewing Horikawa surfaces as symplectic 4--manifolds it is natural (and
desirable) to consider the individual Dehn twists arising in the monodromy,
even though the critical values of the pencils are not
pairwise distinct.
\end{remark}

We now consider the fiber of $f$ through the point $p_j$.
Denote by $\gamma_j$ the loop around $f(p_j)$ chosen as part of our
fixed collection of generators of $\pi_1(\DD\setminus\crit(\tilde{f}))$,
and by $\DD_j$ the disc bounded by $\delta_j$. Since $f$ has no critical
values in $\DD_j$, we can
identify $f^{-1}(\DD_j)$ with $\DD_j\times \Sigma$.
The branch curve $C$ intersects the fiber of $f$ above any point
$z\in\DD_j\setminus \{f(p_j)\}$
transversely in $d$ distinct points, which determines an element $\sigma(z)$
of the configuration space $\mathcal{C}_d(\Sigma)$ (using the trivialization
of $f$). 
The manner in which two of the points in these configurations
converge to each other (while remaining distinct from the others) as $z$
approaches $f(p_j)$ is encoded by a {\it vanishing arc} $\eta_j\subset
\Sigma\setminus \{x_1^0,\dots,x_d^0\}$, with end points in $\{x_1^0,\dots,x_d^0\}$
(recall that we denote by $\{x_1^0,\dots,x_d^0\}\subset\Sigma$
the configuration above the base point). 
Local models for $C$ and $f$ near $p_j$ are given by the plane curve 
$\{y^2=x\}\subset \mathbb{C}^2$ and the projection to the first coordinate.
Therefore, as one moves around $f(p_j)$, the two end
points of the vanishing arc
are exchanged by a counterclockwise half-twist. It follows that the {\it braid
monodromy}\/ of $C$ along $\gamma_j$, ie, the element of $B_d(\Sigma)$
determined by the configurations $\{\sigma(z),\ z\in \gamma_j\}$, is precisely
the half-twist along the vanishing arc $\eta_j$.
By construction, the monodromy of $\tilde{f}$ along $\gamma_j$ is the image
of this half-twist under the lifting homomorphism,
ie, the Dehn twist along $\pi^{-1}(\eta_j)$.

In conclusion, we have proved:

\begin{proposition}\label{prop:lifting}
The vanishing cycles of $\tilde{f}$ are exactly the preimages by
$\pi\co \tilde{\Sigma} \to\Sigma$ of
$(a)$ the vanishing cycles of $f$, and 
$(b)$ the vanishing arcs of the branch curve $C$.
\end{proposition}

\subsection{Horikawa surfaces as genus 2 fibrations}\label{ss:horig2}

The composition of the covering map $\pi_1\co X_1\to\CP^1\times\CP^1$
with the projection to the first factor defines a genus 2 fibration
$\varphi_1\co X_1\to \CP^1$. A generic choice of the branch curve $C_1\subset\CP^1
\times\CP^1$ ensures that every fiber of the projection to the
first factor is tangent to $C_1$ in at most one point, and that every such
tangency is nondegenerate. We can then derive the monodromy of the Lefschetz
fibration $\varphi_1$ from the braid monodromy of the curve $C_1$, as in
\fullref{ss:lifting}. It is easy to check (see eg \cite[Section 3]{Ag2})
that the braid monodromy of $C_1$ corresponds to the
factorization $$(\sigma_1\cdot\sigma_2\cdot \sigma_3\cdot \sigma_4\cdot
\sigma_5\cdot \sigma_5\cdot \sigma_4\cdot \sigma_3\cdot \sigma_2\cdot
\sigma_1)^{12}=1$$
in the spherical braid group $B_6(S^2)$, where $\sigma_1,\dots,\sigma_5$ are
the standard Artin generators (half-twists exchanging two consecutive
points).

\begin{figure}[ht!]
\centering
\includegraphics{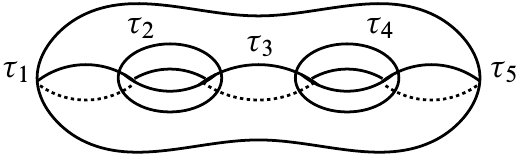}
\caption{Standard generators of $\Map_2$}\label{fig:genus2}
\end{figure}

Similarly, the composition of the covering map $\pi_2\co X_2\to \FF_6$ with
the projection $pr\co \FF_6\to \CP^1$ defines a genus 2 fibration
$\varphi_2\co X_2\to\CP^1$; the braid monodromy of the branch curve
$\Delta_\infty\cup C_2$ can again be expressed in terms of a factorization
in $B_6(S^2)$, namely
$$(\sigma_1\cdot\sigma_2\cdot\sigma_3\cdot\sigma_4)^{30}=1.$$
Using the fact that the half-twists $\sigma_1,\dots,\sigma_5$ lift to the
Dehn twists  $\tau_1,\dots,\tau_5\in \Map_2$ represented in \fullref{fig:genus2} (the standard generators of $\Map_2$), we obtain formulas
for the monodromies of the Lefschetz fibrations $\varphi_1$ and $\varphi_2$.
For another derivation of these formulas, see the work of Fuller \cite{Fu};
see also \cite[Section 4]{ST1} for related considerations.

\begin{proposition}[Fuller]\label{prop:fuller}
$X_1$ and $X_2$ admit genus $2$ Lefschetz fibrations with $120$ singular
fibers; the corresponding monodromy factorizations in $\Map_2$ are
$(\tau_1\cdot\tau_2\cdot\tau_3\cdot\tau_4\cdot\tau_5\cdot
\tau_5\cdot\tau_4\cdot\tau_3\cdot\tau_2\cdot\tau_1)^{12}=1$ and
$(\tau_1\cdot\tau_2\cdot\tau_3\cdot\tau_4)^{30}=1$, respectively.
\end{proposition}

It is easy to see that the Lefschetz fibrations $\varphi_1$ and $\varphi_2$ are
not isomorphic. For example, their monodromy groups are different: the
monodromy of $\varphi_1$ surjects onto $\Map_2$, while that of
$\varphi_2$ takes values in the subgroup of $\Map_2$ generated by 
$\tau_1$, $\tau_2$, $\tau_3$ and $\tau_4$ (this is a proper subgroup
since its image under the natural surjective homomorphism $\Map_2\to 
\mathfrak{S}_6$ mapping $\tau_i$ to $(i,i\!+\!1)$ is $\mathfrak{S}_5
\subsetneq \mathfrak{S}_6$). However, this difference sheds very little
light on the structure of $X_1$ and $X_2$ as symplectic
4--manifolds; in fact, it can be understood in terms of elementary
topological considerations. We start with a remark.

\begin{remark}\label{rmk:blowup}
As a smooth 4--manifold, $X_1$ can also be constructed as follows:
in $\CP^1\times\CP^1$, consider a configuration of 6 ``horizontal'' lines
$H_i=\CP^1\times \{b_i\}$ and 12 ``vertical'' lines $F_j=\{a_j\}\times \CP^1$;
blow up $\CP^1\times\CP^1$ at the 72 intersection points $(a_j,b_i)$, and
denote by $\hat{H}_i$ and $\hat{F}_j$ the
proper transforms of the lines $H_i$ and $F_j$. Then $X_1$ is diffeomorphic
to the double cover of the blowup of $\CP^1\times\CP^1$ branched along 
$\bigcup \hat{H}_i\cup\bigcup \hat{F}_j$. 

Indeed, this follows from
simultaneous resolution of singularities: the double cover of $\CP^1\times\CP^1$
branched along the nodal configuration $\bigcup H_i\cup\bigcup F_j$ is a
singular surface with 72 ordinary double points. The double points can be
either smoothed, which amounts to smoothing of the branch curve in
$\CP^1\times\CP^1$, or blown up, which amounts to blowing up
$\CP^1\times\CP^1$ and taking the proper transform of the branch curve.
Even though these two constructions differ from a symplectic point of view
(blowing up creates symplectic $-2$--spheres, while smoothing creates
Lagrangian $-2$--spheres), the resulting 4--manifolds are diffeomorphic.

The same argument yields an alternative construction of $X_2$ as a double
cover of a blow-up of $\FF_6$.
\end{remark}

The cohomology groups $H^2(X_1,\Z)$ and $H^2(X_2,\Z)$ contain 
rank 2 sublattices $\Lambda_1=\pi_1^*H^2(\CP^1\times\CP^1,\Z)$ and 
$\Lambda_2=\pi_2^*H^2(\FF_6,\Z)$. Even though the lattices 
$H^2(X_1,\Z)$ and $H^2(X_2,\Z)$ (equipped with the intersection pairings)
are isomorphic, and the sublattices $\Lambda_1$ and $\Lambda_2$ (equipped
with the restrictions of the intersection pairings) are also isomorphic,
we claim that the pairs $(H^2(X_i,\Z), \Lambda_i)$ are not
isomorphic. This implies:

\begin{proposition} \label{prop:lattices}
There is no homeomorphism $h\co X_1\to X_2$ such that $h^*(\Lambda_2)=
\Lambda_1$.
\end{proposition}

\proof
Take an element of $\Lambda_1$ of the form $\pi_1^*(p,q)$, where we
implicitly identify $H^2(\CP^1\times\CP^1,\Z)$ with $\Z^2$. Recall
the description of $X_1$ given in \fullref{rmk:blowup}, and consider
the homology classes $A,B\in H_2(X_1,\Z)$ represented by the preimages
of $\hat{H}_1$ and $\hat{F}_1$. By construction, $\langle \pi_1^*(p,q),A
\rangle=p$ and $\langle \pi_1^*(p,q),B\rangle=q$. This implies that, if
$\pi_1^*(p,q)$ is divisible by 2 in $H^2(X_1,\Z)$, then $p$ and $q$ are
both even, and hence $\pi_1^*(p,q)$ is also divisible by 2 in $\Lambda_1$.

On the other hand, let $\alpha\in H^2(\FF_6,\Z)$ be the class Poincar\'e
dual to the exceptional section $\Delta_\infty$, and consider $\pi_2^*\alpha
\in\Lambda_2$: for every class $[C]\in H_2(X_2,\Z)$, we have $$\langle 
\pi_2^*\alpha,[C]\rangle=\langle \alpha,(\pi_2)_*[C]\rangle=[\Delta_\infty]
\cdot (\pi_2)_*[C]=2\,[\pi_2^{-1}(\Delta_\infty)]\cdot [C]$$ (in the last
equality we have used the fact
that $\Delta_\infty$ is a component of the branch curve of $\pi_2$).
This implies that $\pi_2^*\alpha$ is divisible by 2 in $H^2(X_2,\Z)$;
however, $\alpha$ is primitive in $H^2(\FF_6,\Z)$, so $\pi_2^*\alpha$ is
not divisible by 2 in $\Lambda_2$. This completes the proof.
\endproof

The fibers of $\varphi_1$ represent the class $\pi_1^*(0,1)$, while the
fibers of $\varphi_2$ represent the class $\pi_2^*[F]$ where $[F]$ is the
class of the fiber of $\FF_6$. Moreover, the canonical classes of $X_1$
and $X_2$ are $c_1(K_{X_1})=\pi_1^*(1,4)$ and 
$c_1(K_{X_2})=\pi_2^*([\Delta_0]+[F])$.
(Here we are implicitly using the isomorphism between homology and
cohomology given by Poincar\'e duality).

If the two Lefschetz fibrations $\varphi_1$ and
$\varphi_2$ were isomorphic then we would have a diffeomorphism
$h\co X_1\to X_2$ taking the fiber class to the fiber class. Moreover, $h$
would also map the canonical class to the canonical class (this follows
eg from Seiberg--Witten theory, since $\pm K_{X_i}$ are the only basic
classes, and evaluation on the fiber classes shows that the signs are
preserved). Since the fiber classes and canonical classes generate
$\Lambda_1$ and $\Lambda_2$, this contradicts \fullref{prop:lattices}.

\begin{remark}\label{rmk:lag}
These simple topological considerations are at the heart of the problem.
Indeed, it is easy to distinguish $X_1$ and $X_2$ as complex surfaces
because the projections $\pi_1,\pi_2$ and the lattices $\Lambda_1,\Lambda_2$
are naturally determined by the complex geometry of the Horikawa surfaces:
they can, for example, be interpreted in terms of the canonical linear
systems, or in terms of algebraic vanishing cycles for nodal degenerations.
If there were a purely symplectic construction allowing us to
characterize the lattices $\Lambda_1$ and $\Lambda_2$ in terms of
the symplectic topology of $X_1$ and $X_2$ (without using the extra data
provided by the coverings $\pi_i$), then it would follow that $X_1$ and
$X_2$ are different symplectic 4--manifolds. Donaldson has suggested
that one should compare the sets of homology classes realized by embedded
Lagrangian spheres in $X_1$ and $X_2$; it is conjectured that these coincide
with algebraic vanishing cycles, and hence span precisely the orthogonal
complements to $\Lambda_1$ and $\Lambda_2$ (see \fullref{ss:matching} for more
on this topic). However, to this date little progress has been made in this
direction.
\end{remark}

\subsection{Luttinger surgery and partial twistings}\label{ss:luttinger}

We now discuss some properties of braiding constructions and Luttinger surgery
in the context of Lefschetz pencils and double covers. The reader is referred
to \cite{ADK} for more background on these topics (see also \cite{FS}).

\begin{figure}[ht!]
\labellist\small
\pinlabel $A$ [l] at 64 56
\pinlabel $C$ at 107 29
\pinlabel $C'$ at 263 29
\endlabellist
\centering
\includegraphics[scale=0.9]{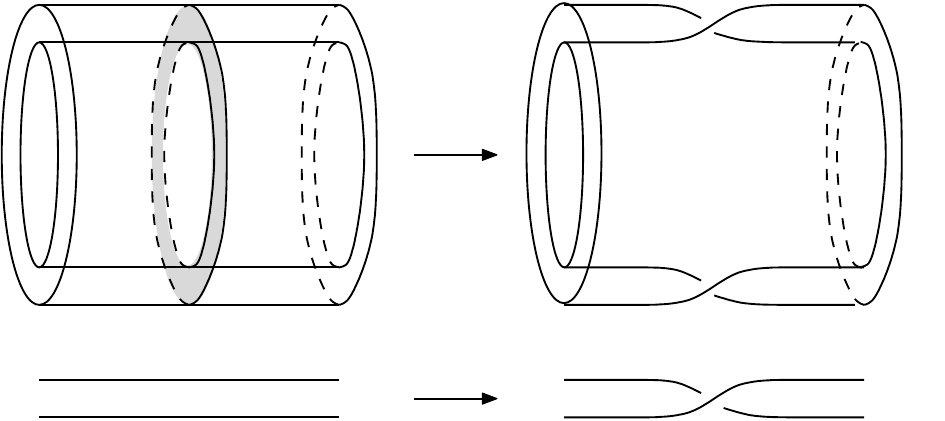}
\caption{Braiding a branch curve}\label{fig:braiding}
\end{figure}

Consider a 2:1 covering $\pi\co X\to Y$ of symplectic 4--manifolds,
branched along a smooth symplectic curve $C\subset Y$. Assume that
we are given a Lagrangian annulus $A$ with interior in $Y\setminus C$
and boundary contained in $C$. Then we can obtain a new symplectic curve
$C'\subset Y$ by {\it braiding} the curve $C$ along the annulus $A$,
in the manner depicted on \fullref{fig:braiding}. Namely, we cut out
a neighborhood $U$ of $A$, and glue it back via a non-trivial diffeomorphism
which interchanges two of the connected components of $C\cap \partial U$,
in such a way that the product of $S^1$ with the trivial two-strand braid
is replaced by the product of $S^1$ with a half-twist (see \cite{ADK} for details).

Braiding the curve $C$ along the Lagrangian annulus $A$ affects the branched
cover $X$ by a {\it Luttinger surgery}\/ along the smooth embedded Lagrangian
torus $T=\pi^{-1}(A)$ \cite{ADK}. This operation consists of cutting out
from $X$ a tubular neighborhood of $T$, foliated by parallel
Lagrangian tori, and gluing it back via a symplectomorphism wrapping the
meridian around the torus (in the direction of the preimage of an arc
joining the two boundaries of $A$), while the longitudes are not affected.

Here we are specifically interested in the situation already considered in
\fullref{ss:lifting}, namely when we are given a Lefschetz pencil (or
fibration) $f\co Y\setminus B\to \CP^1$ with respect to which $C$ lies in a
generic position. Assume that we have a loop $\gamma\in \pi_1(\CP^1\setminus
\crit(f))$ along which the monodromy of $f$ is trivial, and that
$C$ is transverse to the fibers of $f$ in a neighborhood of
$f^{-1}(\gamma)$. Also assume for now that the braid monodromy of $C$ along
$\gamma$ is trivial. 

Consider an arc $\eta$ in the fiber $\Sigma$ of $f$
above a point of $\gamma$, with end points $x',x''\in\Sigma\cap C$.
We can locally identify $Y$ with a product $V_1\times V_2$, where $V_1$
is a neighborhood of $\gamma$ in $\CP^1\setminus \crit(f)$ (such that
$C$ is transverse to the fibers of $f$ over $V_1$), and $V_2$ is
a neighborhood of $\eta$ in $\Sigma$. In this local model, we can assume
that $f$ is the projection to the first factor $V_1$, and that $C$ is the
subset $V_1\times \{x',x''\}\subset V_1\times V_2$. Consider
the annulus $A=\gamma\times \eta$ with boundary in $C$. Standard results
about symplectic structures on Lefschetz fibrations (see eg \cite{Go,GS})
imply that, up to a small perturbation of the symplectic form, we can assume
the annulus $A$ to be Lagrangian. We can then braid $C$ along $A$ as
described above, to obtain a new symplectic curve $C'\subset Y$ which
coincides with $C$ outside of $V_1\times V_2$ and is transverse to the
fibers of $f$ inside $V_1\times V_2$. (In fact, the construction can often
be carried out without perturbing the symplectic form, provided the annulus
$A$ is ``thin'' enough and the neighborhood $V_1$ can be chosen large enough
to allow
$C$ to be distorted in the direction of $\eta$ while remaining symplectic).

The braid monodromy of $C'$ differs from that of $C$ by a {\it partial
conjugation} operation. Namely, let $\DD_\pm$ be the two
components of $\CP^1\setminus V_1$, and choose the base point in
$\CP^1\setminus \crit(f)$ to lie on the boundary of $\DD_+$.
Then $C$ and $C'$ have the same braid monodromy along any loop in
$\DD_+$, but their braid monodromies around points of $\DD_-$ differ
by conjugation by the half-twist $\sigma_\eta$ along $\eta$. In other terms, the
vanishing arcs corresponding to the vertical tangencies of $C$ inside
$f^{-1}(\DD_-)$ are replaced by their images under $\sigma_\eta$.

Since we have only used the local structure near the annulus $A$,
we can in fact relax our assumption concerning the braid monodromy of $C$
along $\gamma$: it is sufficient to assume that this braid monodromy fixes
the arc $\eta$ (ie, that it can be realized by an isotopy of
$\Sigma$ supported in the complement of $V_2$. Similarly, we do not have
to assume that the monodromy of $f$ along $\gamma$ is trivial, we only need
to assume triviality in the considered portion of the fiber (however,
we will only consider situations in which the monodromy of $f$ along
$\gamma$ is trivial away from the boundary of $\Sigma$).

Now consider the double covers $\pi\co X\to Y$ and $\pi'\co X'\to Y$ branched
along $C$ and $C'$, and the Lefschetz pencils $\tilde{f}=f\circ \pi$ and
$\tilde{f}'=f\circ \pi'$. By the result of \cite{ADK}, these differ by a
Luttinger surgery along the
Lagrangian torus $T=\pi^{-1}(A)$, performed in the direction of the
loop $\delta=\pi^{-1}(\eta)\subset\tilde\Sigma$.

The monodromies of $\tilde{f}$ and $\tilde{f}'$
differ by a partial conjugation:

\begin{proposition}\label{prop:partialconj}
The Lefschetz pencils $\tilde{f}$ and $\tilde{f}'$ have the same monodromy
along any loop
in $\DD_+$, but their monodromies around critical values in $\DD_-$
differ by conjugation by the Dehn twist along $\delta$.
\end{proposition}

\proof
The result follows directly from \fullref{prop:lifting}. Indeed,
it is clear that the monodromies over $\DD_+$ coincide. If we consider
a vertical tangency of $C$ in $f^{-1}(\DD_-)$, the braiding operation
replaces the vanishing arc by its image under $\sigma_\eta$; since the
lifting homomorphism maps $\sigma_\eta$ to the Dehn twist $t_\delta$,
the corresponding vanishing cycles of $\tilde{f}$ and $\tilde{f}'$ differ
precisely by $t_\delta$.  Next, consider a critical point of $f$
in $f^{-1}(\DD_-)$. Braiding is only a modification of the curve
$C$, so the vanishing cycles of $f$ are not affected, and neither are the
corresponding vanishing cycles of $\tilde{f}$. 
However, viewing the
half-twist $\sigma_\eta$ as an isotopy of $\Sigma$ supported in $V_2$,
and observing that $[\sigma_\eta]$ is the identity element in the mapping
class group of $\Sigma$,
we can also replace each vanishing cycle of $f$ by its image under
$\sigma_\eta$; after lifting, the corresponding vanishing cycles of
$\tilde{f}'$ become the images of those of $\tilde{f}$ under $t_\delta$
(which, in fact, acts trivially).
\endproof

Another way to interpret these constructions is in terms of {\it twisted
fiber sums}. Namely, $Y$ can be obtained by gluing $Y_+=f^{-1}(\DD_+)$ and
$Y_-=f^{-1}(\DD_-)$ along their boundary via a diffeomorphism
$\psi\co \partial Y_+\to \partial Y_-$ compatible with the Lefschetz fibrations
$f_\pm=f_{|Y_\pm}$; the branch curve
$C$ is also obtained by gluing $C_\pm=C\cap f^{-1}(\DD_\pm)$ along their
boundaries via the diffeomorphism $\psi$. This realizes $(Y,C)$ as
the pairwise fiber sum of $(Y_+,C_+)$ and $(Y_-,C_-)$. If we instead glue
$(Y_+,C_+)$ to $(Y_-,C_-)$ by the diffeomorphism obtained by composing
$\psi$ with the half-twist $\sigma_\eta$ inside each fiber of $f$ above
the boundary (abusing notation we denote this diffeomorphism by
$\sigma_\eta\circ\psi$), we obtain the pair $(Y,C')$, realized as
the twisted fiber sum of $(Y_+,C_+)$ and $(Y_-,C_-)$.

Passing to the double covers, we can view $X$ as a fiber sum
$X_+\cup_{\smash{\tilde\psi}} X_-$, where $X_\pm=\pi^{-1}(Y_\pm)$ and $\tilde\psi$
is a fiber-preserving diffeomorphism which lifts $\psi$. In this language,
the Lefschetz pencil $\tilde{f}'$ on $X'$ is the twisted fiber sum
$X_+\cup_{\smash{\tilde{\psi}'}} X_-$ of $\tilde{f}_+=\tilde{f}_{|X_+}$ and
$\tilde{f}_-=\tilde{f}_{|X_-}$, where
$\tilde{\psi}'=t_\delta\circ \tilde\psi$ is obtained by composing
$\tilde\psi$ with the Dehn twist $t_\delta$ inside each fiber of
$\tilde{f}$. We also say that $\tilde{f}'$ is a ``partial
twisting'' of $\tilde{f}$.

\section{The monodromy of the canonical pencil on $X_2$}\label{s:monx2}

The goal of this section is to compute the monodromy of a canonical pencil
of curves on the Horikawa surface $X_2$ (expressed as a collection of 196
Dehn twists in $\Map_{17,16}$). The reader who does not care about details
of the setup and calculations may skip directly ahead to \fullref{ss:monx2},
where the final formula is given along with the necessary notations.

\subsection{A special configuration}\label{ss:x2notation}

Recall that $X_2$ is the double cover of $\FF_6$ branched along
$\Delta_\infty\cup C_2$, where $\Delta_\infty$ is the exceptional section
and $C_2$ is a smooth curve in the linear system $|5\Delta_0|$. Since
$K_{X_2}=\pi_2^*([\Delta_0]+[F])$, we can obtain a pencil of curves in
the linear system $|K_{X_2}|$ on $X_2$ by taking the preimages of a pencil
of curves in the linear system $|\Delta_0+F|$ on $\FF_6$. 

The connectedness
of the space of generic configurations (which is the complement of a divisor
in some projective variety) implies that the topology of the resulting pencil
of curves on $X_2$ does not depend on the choices made, as long as the curve $C_2$ and the
chosen pencil on $\FF_6$ are in general position with respect to each other.
Accordingly, we will choose a particular configuration for which the
monodromy calculations are manageable.

The Hirzebruch surface $\FF_6=\PP(\OO\oplus \OO(6))$ can be thought of
as a fiberwise compactification
of the line bundle $\OO(6)$ over $\CP^1$; in this sense, $\Delta_0$ is the
zero section, and $\Delta_\infty$ is the section at infinity. Moreover, we
think of $\CP^1$ as $\C\cup\{\infty\}$, and trivialize $\OO(6)$ over $\C$
by means of a holomorphic section vanishing with order 6 at infinity.
In this trivialization, sections of $\OO(6)$ are represented by polynomials
of degree at most $6$ in one complex variable.

We will be considering a pencil of curves in the linear system $|\Delta_0+F|$,
with base points $p_1,\dots,p_8\in\FF_6$; this pencil can be viewed
equivalently as a family of curves $\Sigma_\alpha\subset\FF_6$,
$\alpha\in\CP^1$, or as a map $f\co \FF_6\setminus \{p_1,\dots,p_8\}\to\CP^1$.
We choose the base points $p_1=(z_1,0)$, \dots,
$p_7=(z_7,0)$ on the zero section, and $p_8=(z_8,\epsilon)$ close to the
zero section ($\epsilon$ is a small nonzero constant). 
We will set things up in such a way that all the interesting phenomena
happen in the real part of $\FF_6$ (thus making it easier to visualize the
monodromy). Accordingly, we choose the constants $z_1,\dots,z_8,\epsilon$
to be real numbers, with $z_1<\dots<z_8$ and $\epsilon<0$.

Each curve of the pencil through $p_1,\dots,p_8$ intersects $\Delta_\infty$
in exactly one point, and so we can
parameterize the pencil by its restriction to $\Delta_\infty$: namely,
for each $\alpha\in\CP^1$, we call $\Sigma_\alpha$ the curve in the pencil
which passes through the point $(\alpha,\infty)\in\Delta_\infty$.
If $\Sigma_\alpha$ is smooth, then it can be viewed as the graph of a
meromorphic section $s_\alpha$ of $\OO(6)$ with a simple pole at
$\alpha$ and zeroes at $z_1,\dots,z_7$; in the given trivialization of
$\OO(6)$, this section is given by the formula
$$s_\alpha(z)=\epsilon'\,\frac{z_8-\alpha}{z-\alpha}\,\prod_{i=1}^7(z-z_i),$$
where $\epsilon'=\epsilon/\prod (z_8-z_i)$.
By projecting to $\CP^1$, we can identify each smooth fiber of 
$f$ (ie, $\Sigma_\alpha$
with the base points removed) with $\CP^1\setminus \{z_1,\dots,z_8\}$.
Accordingly, the monodromy of $f$ takes values in the
mapping class group $\Map_{0,8}$ of $\CP^1$ with small discs around each
$z_i$ removed.

It is easy to check that $f$ has 8 singular fibers, corresponding to the values
$\alpha=z_1,\dots,z_8$; for $\alpha=z_i$, the curve $\Sigma_\alpha$ consists
of two components: the fiber of $\FF_6$ over $z_i$, and the unique curve in
$|\Delta_0|$ passing through all $p_j$, $j\neq i$ (for $i=8$ this is the
zero section). The base point $p_i$ lies in the fiber component, and the
seven other base points lie in the section component; hence, the vanishing
cycle of $f$ at $\alpha=z_i$ is a boundary curve (a circle separating $z_i$
from the other punctures).

Choosing $\epsilon$ small enough ensures that, outside of a fixed neighborhood
of $\Delta_0$, the pencil of curves $(\Sigma_\alpha)_{\alpha\in\CP^1}$ is
arbitrarily close to the standard fibration $\FF_6\to\CP^1$. In fact,
each $\Sigma_\alpha$ lies in a small neighborhood of $\Delta_0\cup (\{\alpha\}
\times\CP^1)$; therefore, by choosing the curve $C_2$ (the ``main'' component of
the branch curve of $\pi_2$) transverse to the zero section, we can obtain
an explicit description of its behavior with respect to the pencil $f$.

Choose real numbers $q_1,\dots,q_6$ and $r_1,\dots,r_6$ such that
$$q_1<z_1<r_1<q_2<z_2<r_2<\dots<q_6<z_6<r_6<z_7<z_8,$$
and consider the graph of the holomorphic section
$u(z)=\prod_{i=1}^6 (z-q_i)$ of $\OO(6)$.
We can construct nearby curves in the real part of the pencil of curves in
$|\Delta_0|$, passing through the six points $(r_i,u(r_i))$, by
considering the graphs of holomorphic sections of $\OO(6)$ of the form
$$u(z)+\lambda\,\prod_{i=1}^6 (z-r_i),$$ with $\lambda\in\R$ small.
Consider five such
holomorphic sections $u_1,\dots,u_5$ of $\OO(6)$: their graphs
$\Gamma_1,\dots,\Gamma_5$ intersect each other transversely at
$(r_i,u(r_i))$ ($1\le i\le 6$), and intersect the zero
section transversely at real points near $q_1,\dots,q_6$. We define
$\Gamma_j\cap \Delta_0=\{q_{1,j},\dots,q_{6,j}\}$, with $q_{i,j}\approx q_i$.
It is easy to check that the various points $q_{i,j}$ are in the same
order near each $q_i$, so we can assume that $q_{i,1}<\dots<q_{i,5}$ for all
$i$. Choosing the perturbations small enough, we can also assume that
$q_{i,1}>r_{i-1}$ and $q_{i,5}<z_i$.

\begin{figure}[ht!]
\centering
\includegraphics{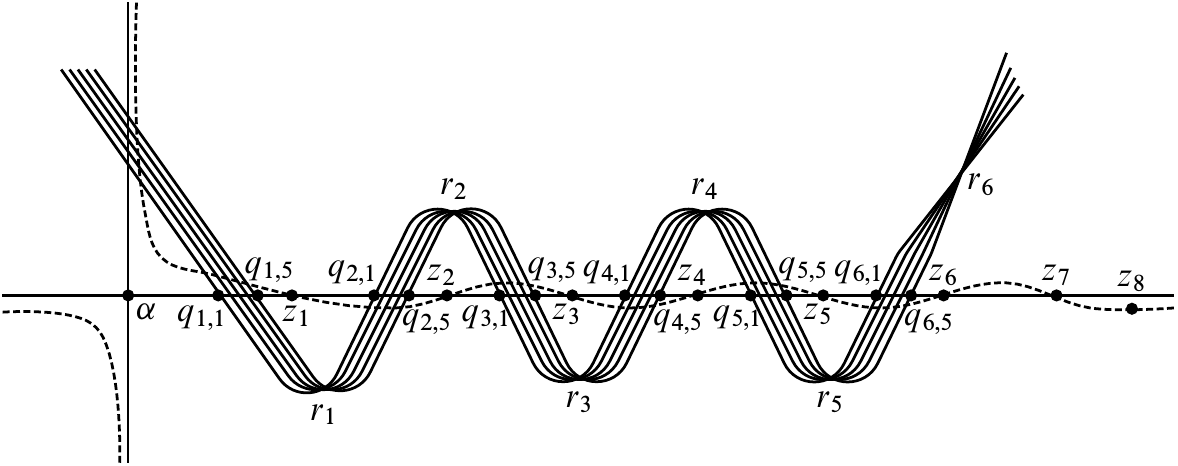}
\caption{The curves $C_2$ (solid) and $\Sigma_\alpha$ (dashed) in $\FF_6$}\label{fig:x2setup}
\end{figure}

The configuration $\Gamma_1\cup\dots\cup \Gamma_5$ can be smoothed to a
nearby curve in the linear system $|5\Delta_0|$, which we take as our choice
for $C_2$. Since this smoothing
can be realized by an arbitrarily small perturbation, we can ensure that
$C_2$ is contained in an arbitrarily small neighborhood of
$\bigcup \Gamma_i$, and arbitrarily $C^1$--close to $\bigcup \Gamma_i$ outside
of an arbitrarily small neighborhood of the points $(r_i,u(r_i))$.
Moreover, we choose the smoothing perturbation to be real and generic,
so that $C_2$ is defined by an equation with real coefficients, and its tangencies
with the fibers of $\FF_6$ are nondegenerate, real, and lie in distinct fibers near
$r_i$. Finally, the points of $C_2\cap \Delta_0$ lie arbitrarily close to
those of $\Gamma_j\cap \Delta_0$; changing our notation, we will again call
them $q_{i,j}$, $1\le i\le 6$, $1\le j\le 5$, and observe that we still
have $r_{i-1}<q_{i,1}<\dots<q_{i,5}<z_i$.
Our choices for $C_2$ and $\Sigma_\alpha$ are summarized on \fullref{fig:x2setup}.

\subsection{More notations and conventions}
We need to study the braid monodromy of the branch curve $\Delta_\infty\cup C_2$
with respect to the pencil $f$, ie, the manner in which the intersection
points of $\Delta_\infty\cup C_2$ with $\Sigma_\alpha$ depend on the choice
of $\alpha$. Recall that we can trivialize the pencil $f$ (except at its
singular fibers) by using the standard projection $pr\co \FF_6\to\CP^1$,
which allows us to identify the complement of the base points
in $\Sigma_\alpha$ with $\CP^1\setminus \{z_1,\dots,z_8\}$. 

Assume that $\alpha$ is not too close to any of the
special values $z_i,q_i,r_i$. Then the projections to $\CP^1\setminus
\{z_1,\dots,z_8\}$ of the 36 points where $\Sigma_\alpha$
intersects $\Delta_\infty\cup C_2$ all lie in a small neighborhood of 
$\{\alpha,q_1,\dots,q_6\}$, and can be labelled in a simple manner
according to their respective positions:
\begin{itemize}
\item the intersection between $\Sigma_\alpha$ and $\Delta_\infty$
takes place at $\alpha$.
\item the 5 intersections between $C_2$ and the ``vertical'' part of
$\Sigma_\alpha$ take place at $\tilde{\alpha}_1,\dots,\tilde{\alpha}_5$
in a neighborhood of $\alpha$. We label them in such a way that, upon
deforming $C_2$ back to the nearby singular configuration $\Gamma_1\cup\dots\cup
\Gamma_5$, $\tilde{\alpha}_i$ corresponds to an intersection of
$\Sigma_\alpha$ with~$\Gamma_i$.
\item the 30 intersections between $C_2$ and the ``horizontal'' part of
$\Sigma_\alpha$ take place at $\tilde{q}_{1,1},\dots,\tilde{q}_{6,5}$, where
each $\tilde{q}_{i,j}$ is close to $q_{i,j}$.
\end{itemize}

In order to define braid monodromy and vanishing arcs, we need to fix a
reference fiber of the pencil $f$, ie, some reference value $\alpha_0$,
and arcs from this reference value to the various values of $\alpha$ for
which $\Delta_\infty\cup C_2$ is tangent to $\Sigma_\alpha$. We choose
the reference value $\alpha_0$ to be a sufficiently negative real number, so
that $\alpha_0\ll q_{1,1}$. The respective positions of $\Sigma_{\alpha_0}$ and
$\Delta_\infty\cup C_2$ are then as pictured on \fullref{fig:x2setup}.
In particular, the images under the projection from $\Sigma_{\alpha_0}$
to $\CP^1$ of the 36 intersection points and the 8 base points are all
real, and in the order
\begin{multline*}
\alpha<\tilde\alpha_5<\dots<\tilde\alpha_1<\tilde{q}_{1,1}<\dots<
\tilde{q}_{1,5}<z_1<\tilde{q}_{2,1}<\dots<\tilde{q}_{2,5}<z_2<\dots\\
\dots<z_5<\tilde{q}_{6,1}<\dots<\tilde{q}_{6,5}<z_6<z_7<z_8.
\end{multline*}
To determine the braid monodromy of $\Delta_\infty\cup C_2$, we consider
what happens to these various intersections as the value of $\alpha$
increases along the real axis from $\alpha_0$ to a large positive value.
As we will see below, there are in total 180 values of $\alpha$ for which
the curve $\Sigma_\alpha$ is tangent to $C_2$, in addition to the 8 values
of $\alpha$ for which $\Sigma_\alpha$ is nodal. Our convention will be that
we determine the monodromy around each critical value $\alpha_{cr,i}$ by
considering a loop in $\pi_1(\C\setminus \{\alpha_{cr,j}\}, \alpha_0)$ 
constructed as follows: choose a point $\alpha'_i$ on the real axis just
to the left of the critical value $\alpha_{cr,i}$, and an arc $\eta_i$
joining $\alpha_0$ to $\alpha'_i$ inside the upper half-plane (ie, passing
{\it above} all the critical values between $\alpha_0$ and $\alpha_{cr,i}$);
then we consider the loop $\gamma_i$ obtained by composing the arc $\eta_i$ from
$\alpha_0$ to $\alpha'_i$, a small circle around $\alpha_{cr,i}$
(counterclockwise), and the arc $\eta_i^{-1}$ back to $\alpha_0$.
This choice ensures in particular that, if we order the critical values
$\alpha_{cr,i}$ in increasing order along the real axis, the loops $\gamma_i$
form an ordered collection of generators for $\pi_1(\C\setminus
\{\alpha_{cr,j}\},\alpha_0)$.

Our calculation of braid monodromy relies on the following ingredients:

\begin{enumerate}
\item
The configuration of intersection points for $\alpha=\alpha'_i$
determines readily the vanishing arc at $\alpha_{cr,i}$: namely, in the
nearby fiber $\Sigma_{\alpha'_i}$ the vanishing arc is simply a straight
line segment joining the two intersection
points which approach each other as $\alpha\to\alpha_{cr,i}$.
\item
The vanishing arc in the reference fiber $\Sigma_{\alpha_0}$
is obtained from the local configuration in $\Sigma_{\alpha'_i}$ by
transporting along the arc $\eta_i^{-1}$, or equivalently, along a
succession of counterclockwise half-circles around all the critical values
$\alpha_{cr,j}<\alpha_{cr,i}$. (As a general principle, the main feature
is that the intersection points labelled
$\alpha,\tilde\alpha_5,\dots,\tilde\alpha_1$ are moved counterclockwise
back to the leftmost positions, since these intersection points stay
close to $\alpha$ while the others remain close to $q_1,\dots,q_5$).
\item
The configuration of intersection points for a value of
$\alpha$ on the real axis just to the right of $\alpha_{cr,i}$
can be deduced from the configuration at $\alpha=\alpha'_i$ by applying
the ``square root'' of the monodromy clockwise around $\alpha_{cr,i}$, namely a
clockwise 90 degree rotation of the two end points of the local
vanishing arc.
\end{enumerate}

\subsection{The braid monodromy}\label{ss:calcx2}

\begin{figure}[ht!]
\centering
\includegraphics{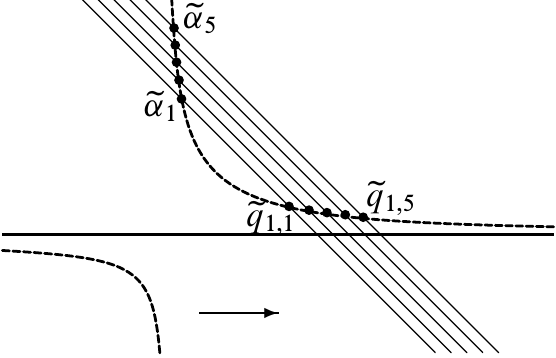}
\caption{The configuration near $q_1$}\label{fig:nearq1}
\end{figure}

Upon increasing $\alpha$ along the real axis, the first critical values
encountered lie near $q_1$. The braid monodromy near $q_1$ can be understood
by following the approach outlined above. The real part of the local
configuration looks as in \fullref{fig:nearq1} (a translating
hyperbola passing through five parallel lines); in particular, near $q_1$
there are 10 values $\alpha_{cr,1}<\dots<\alpha_{cr,10}$ of $\alpha$ for
which $\Sigma_\alpha$ is tangent to $C_2$. 

\begin{figure}[ht!]
\centering
\includegraphics{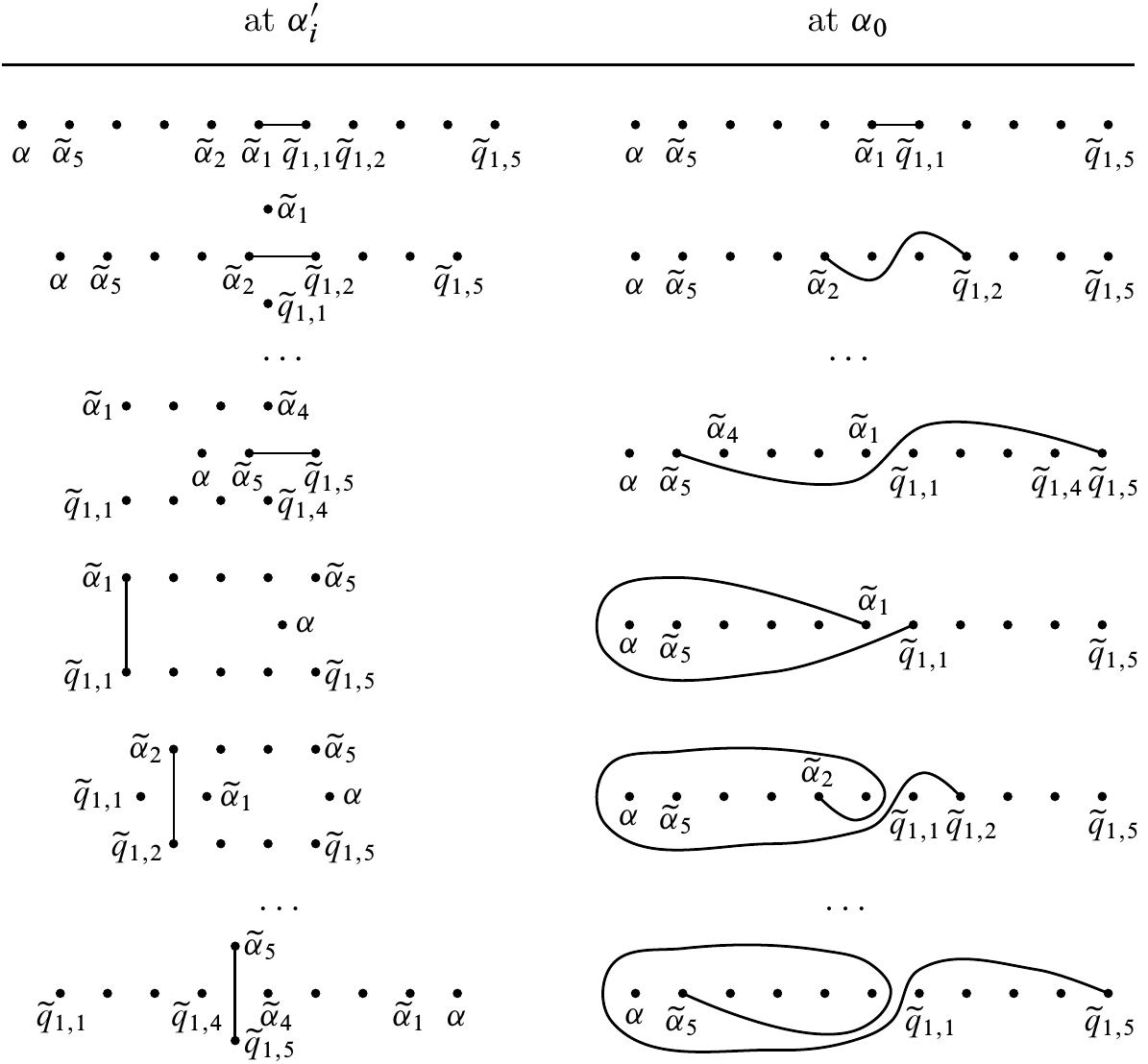}
\caption{The vanishing arcs near $q_1$}\label{fig:monoq1}
\end{figure}

The local configurations of the intersection points in the fibers
immediately to the left of each critical value $\alpha_{cr,1},\dots,
\alpha_{cr,10}$ are shown in \fullref{fig:monoq1} (left), together
with the corresponding vanishing arcs. Transporting these vanishing arcs
back to the reference fiber $\Sigma_{\alpha_0}$ (by going counterclockwise
around the previous critical values) yields the vanishing arcs represented
in \fullref{fig:monoq1} (right), which determine the braid monodromy of
$\Delta_\infty\cup C_2$ near~$q_1$.

The product of the 10 half-twists along these vanishing arcs is the braid
which translates the disc $\DD$ containing $\alpha,\tilde\alpha_5,\dots,
\tilde\alpha_1$ counterclockwise around the points
$\tilde{q}_{1,1},\dots, \tilde{q}_{1,5}$ by 360 degrees,
while simultaneously rotating the interior of $\DD$ clockwise by 360 degrees,
and rotating the interior of the smaller disc $\DD'$ containing
$\tilde\alpha_5,\dots,\tilde\alpha_1$ clockwise by another 360 degrees.
This can be checked either by direct calculation, or more
geometrically by analyzing the behavior of $\alpha,\tilde\alpha_5,\dots,
\tilde\alpha_1$ as $\alpha$ moves along the boundary of a larger circle
enclosing $q_{1,1},\dots,q_{1,5}$ (it is clear that $\DD$ is moved
counterclockwise around $\tilde{q}_{1,1},\dots,\tilde{q}_{1,5}$;
a careful analysis of
the respective positions of the $\tilde\alpha_i$ relatively to $\alpha$ 
shows that the motion within $\DD$ is as claimed).

\begin{figure}[ht!]
\centering
\includegraphics{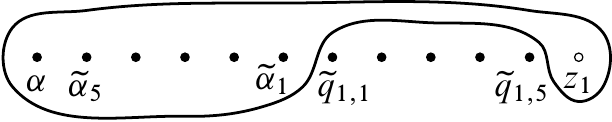}
\caption{The vanishing cycle at $z_1$}\label{fig:monoz1}
\end{figure}

The next contribution to monodromy occurs at $\alpha=z_1$; the curve $\Sigma_\alpha$
is then reducible, with one component (the fiber of $\FF_6$ above $z_1$)
containing the intersection points labelled $\alpha,\tilde{\alpha}_5,\dots,
\tilde{\alpha}_1$ and the base point $z_1$, while the other component
contains all the other intersection points and base points. Hence, in a
nearby smooth fiber the vanishing cycle is a simple closed curve
around the points
$\tilde\alpha_5,\dots,\tilde\alpha_1,\alpha,z_1$ (which are adjacent in that order
along the real axis). Transporting things back to the
reference fiber $\Sigma_{\alpha_0}$ in the prescribed manner, using the
square root of the monodromy around $q_1$ determined above, yields the
vanishing cycle represented in \fullref{fig:monoz1}.

Next, we consider the monodromy near $\alpha=r_1$. Because the point
$(r_1,u(r_1))$
lies away from the zero section, the tangencies between $C_2$ and
$\Sigma_\alpha$ occur in the portion of $\Sigma_\alpha$ which lies close
to the fiber $\{\alpha\}\times\CP^1$. In particular, the braid
monodromy consists of half-twists supported in a small disc containing
the points $\tilde\alpha_1,\dots,\tilde\alpha_5$, and can be understood in
terms of a local model in which a moving fiber intersects a
degree 5 curve obtained by smoothing a configuration of five concurrent lines.

\begin{figure}[ht!]
\centering
\includegraphics{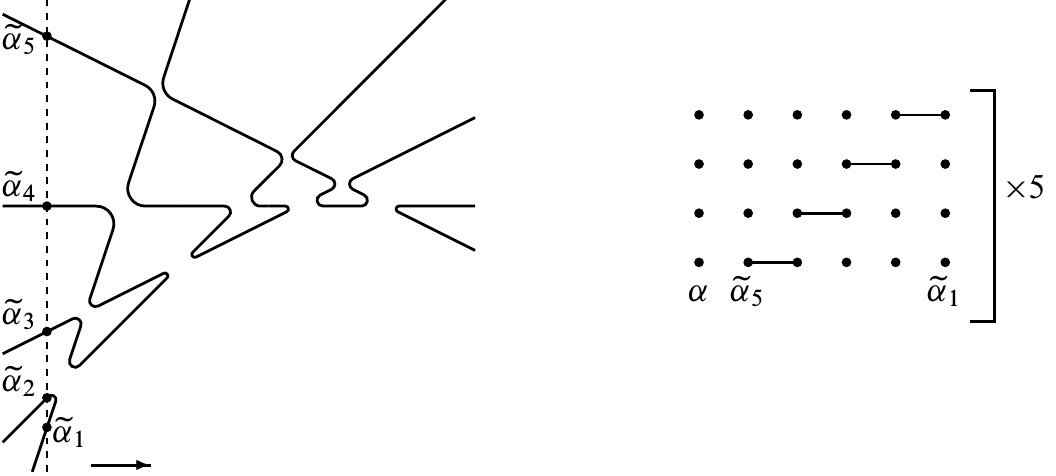}
\caption{The monodromy near $r_1$}\label{fig:monor1}
\end{figure}

Hence, in a neighborhood of $r_1$, the braid monodromy of
$\Delta_\infty\cup C_2$ with respect to the pencil $f$ coincides with
the braid monodromy of a smooth algebraic plane curve of degree 5 with
respect to a generic linear projection. The braid monodromy
of smooth algebraic plane curves has been studied extensively,
and there are various well-known formulas (Hurwitz equivalent
to each other, of course); see eg \cite{Moi}. For completeness,
we outline one possible approach:
first deform the curve so that it lies close to a
configuration of five lines in general position, as in \fullref{fig:monor1}.
Consistently with the choice we have made so far,
consider loops that reach each critical value $\alpha_{cr,i}$ via an
arc in the upper half-plane. In that case, the braid monodromy can be
read off from \fullref{fig:monor1} using the same method as previously;
the braid
monodromy factorization that arises in this way is a product of 20
half-twists,
$$\prod_{i=1}^4 \prod_{j=i+1}^5 (\sigma_{i,j}\cdot\sigma_{i,j}),$$
where $\sigma_{ij}$ is the half-twist along an arc that joins
$\tilde\alpha_i$ to $\tilde\alpha_j$ passing below the real axis
(see \cite{Moi} for a careful derivation of this formula).
However, it is well-known (see eg \cite{Moi}) that this expression
is Hurwitz equivalent to the simpler expression
$$(\sigma_{1,2}\cdot\sigma_{2,3}\cdot\sigma_{3,4}\cdot\sigma_{4,5})^5.$$
In other terms, if we change our choice of ordered collection of generators
for $\pi_1(\C\setminus\{\alpha_{cr,j}\},\alpha_0)$, we can assume that the
20 vanishing arcs near $r_1$ are as pictured in \fullref{fig:monor1}
(right). 

The product of all the monodromies encountered so far (near $q_1$, $z_1$
and $r_1$) is simply the braid which moves the disc $\DD$ containing
$\alpha,\tilde{\alpha}_5,\dots,\tilde{\alpha}_1$ by 360 degrees
counterclockwise around the points $\tilde{q}_{1,1},\dots,\tilde{q}_{1,5},z_1$.
Moreover, for a real value of $\alpha$ such that $r_1\ll \alpha\ll q_2$,
the intersections of $\Sigma_\alpha$ with $\Delta_\infty\cup C_2$ are all
real, and in the order
\begin{multline*}
\tilde{q}_{1,1}<\dots<\tilde{q}_{1,5}<z_1<\alpha<\tilde\alpha_5<
\dots<\tilde\alpha_1<\tilde{q}_{2,1}<\dots<\tilde{q}_{2,5}<z_2<\dots\\
\dots<z_5<\tilde{q}_{6,1}<\dots<\tilde{q}_{6,5}<z_6<z_7<z_8.
\end{multline*}
Hence, the local pictures near $q_2$, $z_2$ and $r_2$ are exactly the
same as near $q_1$, $z_1$ and $r_1$ respectively (except that the
configurations are flipped in the vertical direction, which does not change
anything since we characterize intersection points in terms of their
projections to the horizontal direction). When we transport the local
monodromies back to the reference fiber $\Sigma_{\alpha_0}$ by an arc
that goes counterclockwise around $q_1$, $z_1$ and $r_1$, the points
$\alpha,\tilde{\alpha}_5,\dots,\tilde{\alpha}_1$ are moved back to the
left of $\tilde{q}_{1,1},\dots,\tilde{q}_{1,5},z_1$ by a 180 degree
counterclockwise motion around these points. For example, the vanishing arcs
near $q_2$ look identical to those near $q_1$ except that they connect
$\tilde{\alpha}_j$ to $\tilde{q}_{2,j}$ by passing above the points
$\tilde{q}_{1,1},\dots,\tilde{q}_{1,5},z_1$.

\begin{figure}[ht!]
\centering
\includegraphics{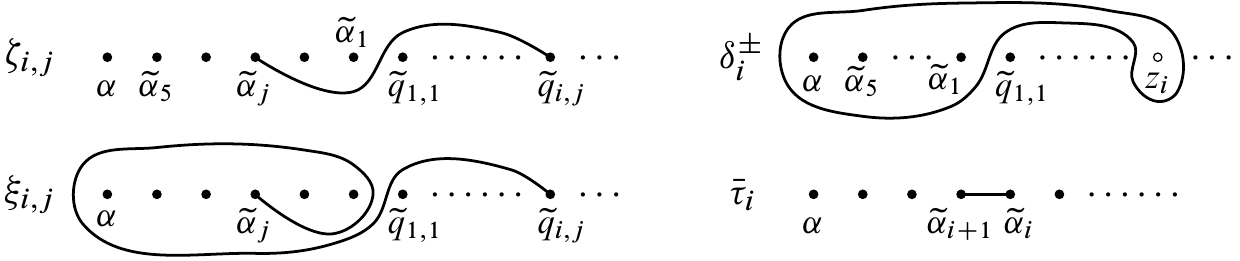}
\caption{The vanishing cycles of the canonical pencil on $X_2$}
\label{fig:monx2}
\end{figure}

The same argument holds for the monodromies near $q_i$, $z_i$ and
$r_i$ for $i\ge 3$; hence, we have now determined 180 vanishing arcs
for $\Delta_\infty\cup C_2$ (10 at each $q_i$ and 20 at each $r_i$),
and 8 vanishing cycles of $f$; see \fullref{fig:monx2}. 
Using \fullref{prop:lifting}, these calculations yield
196 vanishing cycles for the pencil $\tilde{f}_2=f\circ \pi_2$ on $X_2$ (one for
each vanishing arc of $\Delta_\infty \cup C_2$ and two for each vanishing
cycle of~$f$). 

\subsection{The monodromy of the canonical pencil on $X_2$}\label{ss:monx2}
\begin{lemma}
A generic curve in the linear system $|K_{X_2}|$ has genus $17$; a generic
pencil of such curves has $16$ base points, and $196$ nodal singularities.
\end{lemma}

\proof
Since $K_{X_2}\cdot K_{X_2}=2[\Delta_0+F]\cdot[\Delta_0+F]=16$, the adjunction formula
yields that a generic curve in $|K_{X_2}|$ has genus $g=1+K_{X_2}\cdot
K_{X_2}=17$, and two such curves intersect in 16 points. By blowing up the
16 base points of a canonical pencil on $X_2$, we obtain a surface
$\hat{X}_2$ with Euler characteristic $e(\hat{X}_2)=e(X_2)+16=132$. This
surface carries a Lefschetz fibration of genus $g=17$, and the Euler
characteristic is related to the number $N$ of nodal singularities by the
classical formula $e(\hat{X}_2)=4-4g+N$, which implies that $N=196$.
\endproof

Hence, the monodromy of a generic pencil of curves
in the linear system $|K_{X_2}|$ can be expressed in terms of 196 Dehn
twists in the mapping class group $\Map_{17,16}$; in particular, this
confirms that all the vanishing cycles of the pencil $\tilde{f}_2=f\circ \pi_2$ have
been accounted for in the above calculations.

Recall that we view the reference fiber $\tilde{\Sigma}$ of $\tilde{f}_2$ as
a double cover of a sphere with 8 punctures $\Sigma=\CP^1\setminus
\{z_1,\dots,z_8\}$, branched in 36 points. Also
recall that, in $\Sigma$, the punctures $\{z_i,\ 1\le i\le 8\}$
and the branch points $\{\alpha,\,\tilde\alpha_j,\,
\tilde{q}_{i,j},\ 1\le i\le 6,\ 1\le j\le 5\}$ all lie on the real
axis, in the order
\begin{multline*}
\alpha<\tilde\alpha_5<\dots<\tilde\alpha_1<\tilde{q}_{1,1}<\dots<
\tilde{q}_{1,5}<z_1<\tilde{q}_{2,1}<\dots<\tilde{q}_{2,5}<z_2<\dots\\
\dots<z_5<\tilde{q}_{6,1}<\dots<\tilde{q}_{6,5}<z_6<z_7<z_8.
\end{multline*}
With this notation, the calculations in \fullref{ss:calcx2} imply:

\begin{theorem}\label{thm:monx2}
Up to global conjugation and Hurwitz equivalence, the monodromy of a
generic pencil of curves in the linear system $|K_{X_2}|$
is expressed by the factorization of the boundary twist into the
following product of $196$ Dehn twists:
$$
\prod_{i=1}^6\Biggl[\,\prod_{j=1}^5\zeta_{i,j}\cdot\prod_{j=1}^5\xi_{i,j}\cdot
\delta_i^+\cdot\delta_i^-\cdot(\bar\tau_1\cdot\bar\tau_2\cdot\bar\tau_3\cdot\bar\tau_4)^5
\Biggr]\cdot \delta_7^+\cdot\delta_7^-\cdot\delta_8^+\cdot\delta_8^-,
$$
where $\zeta_{i,j}$, $\xi_{i,j}$, $\delta_i^\pm$, and $\bar\tau_i$ are the
Dehn twists along the preimages of the arcs and curves in $\Sigma$
represented in \fullref{fig:monx2}.
\end{theorem}

\begin{remark}
Our convention is to write products of elements in the mapping class group 
in the left-to-right order, consistently with the
standard convention for braid groups. Hence, $\varphi_1\cdot \varphi_2$ is the
mapping class represented by the composition $\varphi_2\circ\varphi_1$.
\end{remark}

\section{A symplectic Lefschetz pencil on $X_1$}\label{s:monx1}

\subsection{Horikawa surfaces and Luttinger surgery}\label{ss:horilu}
The Hirzebruch surface $\FF_6$ can be realized as the
fiber sum of the Hirzebruch surfaces $\FF_2$ and $\FF_4$, in such a way
that $\Delta_\infty\cup C_2$ decomposes into the fiber sum of two curves
$D_2\subset \FF_2$ and
$D_4\subset \FF_4$. Here each $D_k$ ($k\in\{2,4\}$) is the disjoint union of
the exceptional section (of square $-k$) and a smooth curve in the linear
system corresponding to five times a section of square $+k$. Hence,
the genus 2 fibration $\varphi_2\co X_2\to\CP^1$ introduced
in \fullref{ss:horig2} (obtained by composing the double cover 
$\pi_2\co X_2\to\FF_6$ with the standard projection $pr\co \FF_6\to\CP^1$)
is actually the fiber sum of two genus 2 fibrations similarly defined on
the double covers of $\FF_2$ and $\FF_4$ branched along $D_2$ and $D_4$.

With the notations of \fullref{ss:x2notation}, let $\gamma$ be a loop in
$\CP^1\setminus\crit(\varphi_2)$ which bounds a disc
$\DD_-\subset\CP^1$ containing the points $q_1,z_1,r_1,q_2,z_2,r_2$,
while the points $q_i,z_i,r_i$ for $i\ge 3$ lie in
$\DD_+=\CP^1\setminus \DD_-$. The vertical tangencies of
$\Delta_\infty \cup C_2$ lie near the points $r_i$, $1\le i\le 6$,
and each $r_i$ contributes $(\sigma_1\cdot\sigma_2\cdot\sigma_3\cdot
\sigma_4)^5$ to the braid monodromy factorization. Hence, the fiber sum
decomposition described above corresponds exactly to a decomposition of $\CP^1$
into the two discs $\DD_\pm$. The braid monodromy of
$\Delta_\infty\cup C_2$ over $\DD_-$ can be represented by the factorization
$(\sigma_1\cdot\sigma_2\cdot\sigma_3\cdot\sigma_4)^{10}$ in $B_6(S^2)$,
while the
braid monodromy over $\DD_+$ is $(\sigma_1\cdot\sigma_2\cdot\sigma_3
\cdot\sigma_4)^{20}$; these two expressions are precisely the braid
monodromy factorizations of $D_2$ and $D_4$ with respect to the natural
projections. A similar property holds for the monodromies of the
corresponding genus 2 Lefschetz fibrations.

The braid monodromy of $\Delta_\infty\cup C_2$ along
$\gamma$ is trivial, and hence preserves any arc in the fiber of $\FF_6$
with end points on the branch curve.  This allows us to apply the
considerations of \fullref{ss:luttinger}. Namely, choose the
reference fiber of $\FF_6$ to be the fiber above a point of $\gamma$,
and in this fiber
let $\eta$ be the supporting arc of the standard half-twist $\sigma_5\in
B_6(S^2)$, so that one end point of $\eta$ lies on $C_2$ and the other lies
on $\Delta_\infty$; we additionally assume that $\eta$ lies away from the
point where the zero section hits the reference fiber. In fact, any other
arc in the reference fiber with one end point on $C_2$ and the other in
$\Delta_\infty$ and avoiding the zero section would be equally suitable.

As in \fullref{ss:luttinger}, consider the annulus 
$A=\gamma\times\eta\subset \FF_6$ (with one boundary on $\Delta_\infty$
and the other on $C_2$), and its preimage $T=\pi_2^{-1}(A)$, which is a
smoothly embedded torus in $X_2$. Up to suitable perturbations, $A$ and $T$
can be assumed to be Lagrangian.

Let $C'\subset \FF_6$ be the symplectic curve obtained by braiding
$\Delta_\infty\cup C_2$ along the annulus $A$. This curve is connected,
coincides with $\Delta_\infty\cup C_2$ outside of a neighborhood of $A$,
intersects the fibers of $\FF_6$ transversely in a neighborhood of $A$,
and represents the same homology class $5[\Delta_0]+[\Delta_\infty]$ as
$\Delta_\infty\cup C_2$. Recall that we denote by $C_1$ a smooth algebraic
curve of bidegree $(6,12)$ in $\CP^1\times\CP^1$ (the branch curve of 
$\pi_1\co X_1\to\CP^1\times\CP^1$).

\begin{proposition}\label{prop:isotopy}
Equip $\FF_6$ and $\CP^1\times\CP^1$ 
with fixed K\"ahler forms in the cohomology classes Poincar\'e dual to
$[\Delta_0]+[F]$ and $(1,4)$, respectively. Then
there exists a symplectomorphism $\psi\co \CP^1\times\CP^1\to \FF_6$ such
that $\psi(C_1)=C'$.
\end{proposition}

\proof
It is a classical fact that $\FF_6$ and $\CP^1\times\CP^1$ are
diffeomorphic; such a diffeomorphism can be assumed to preserve the
fiber class, and map the homology classes $[\Delta_0]$ and
$[\Delta_\infty]$ to $(1,3)$ and $(1,-3)$ respectively.
Moreover, it is well-known that $\FF_6$ and $\CP^1\times\CP^1$ equipped
with the chosen K\"ahler forms are in fact symplectomorphic, even though
their complex structures are different (see eg \cite[Section 9.4]{McS};
see also \cite{McD}).

Hence, we can think of $C'$ as a connected symplectic curve in
$\CP^1\times\CP^1$, representing the homology class 
$5\cdot(1,3)+(1,-3)=(6,12)$. We now appeal to the following isotopy
result, due to Siebert and Tian \cite{ST}:

\begin{theorem}[Siebert--Tian]\label{thm:ST}
Let $\Sigma\subset \CP^1\times\CP^1$ be a connected symplectic submanifold,
such that the intersection number between $\Sigma$ and the fiber class is
at most $7$. Then $\Sigma$ is symplectically isotopic to a holomorphic
curve.
\end{theorem}

Since any two smooth holomorphic curves in the homology class $(6,12)$ are
mutually isotopic, we conclude that $C'$ is
symplectically isotopic to the branch curve $C_1$ of $\pi_1\co X_1\to\CP^1\times
\CP^1$. This implies the existence of a symplectomorphism of
$\CP^1\times\CP^1$ (in fact, an isotopy, see eg 
\cite[Proposition 0.2]{ST}) which maps $C_1$ to $C'$. 
\endproof

\begin{theorem}\label{thm:horilu}
The manifold obtained from $(X_2,\omega_2)$ by Luttinger surgery along
the torus $T$ in the direction of $\pi_2^{-1}(\eta)$ is symplectomorphic
to $(X_1,\omega_1)$.
\end{theorem}

\proof 
As we have seen in \fullref{ss:luttinger}, the symplectic manifold 
$(X',\omega')$ obtained from $(X_2,\omega_2)$
by Luttinger surgery along $T$ is precisely the double cover of
$\FF_6$ branched along $C'$. In fact, recall that
the double cover of a symplectic 4--manifold branched along a symplectic
curve carries a natural symplectic structure, canonically determined up to
symplectomorphism (see eg \cite[Proposition 3.2]{AS}). The symplectic
forms $\omega'$ and $\omega_1$ are precisely those induced on the
double covers $X'$ and $X_1$ by the chosen K\"ahler forms on
$\FF_6$ and $\CP^1\times\CP^1$. Hence, the symplectomorphism $\psi$ given
by \fullref{prop:isotopy} can be lifted to a
symplectomorphism from $(X',\omega')$ to $(X_1,\omega_1)$.
\endproof

Since the Luttinger surgery operation which yields $X_1$ from $X_2$ is
carried out in a manner compatible with the genus 2 fibration
$\varphi_2\co X_2\to\CP^1$, it yields a symplectic Lefschetz fibration 
$\varphi'\co X'=X_1\to\CP^1$, whose monodromy differs from that
of $\varphi_2$ in the manner described by \fullref{prop:partialconj}.
Together with \fullref{prop:fuller}, this
implies that the monodromy of $\varphi'$ is described by the
factorization
$(\tau_1\cdot\tau_2\cdot\tau_3\cdot \tau')^{10}\cdot
(\tau_1\cdot\tau_2\cdot\tau_3\cdot\tau_4)^{20}=1$ in
$\Map_2$, where $\tau'=\tau_5\tau_4\tau_5^{-1}$.
With some work, one can verify that this factorization is Hurwitz equivalent
to the monodromy of $\varphi_1$, which means that the Lefschetz fibrations
$\varphi'$ and $\varphi_1$ are isomorphic. However, such a result also follows
more directly from the work of Siebert and Tian \cite{ST}.

\subsection{A symplectic Lefschetz pencil on $X_1$}
The monodromy of a canonical pencil on $X_1$ can be determined
directly by the same methods as in \fullref{s:monx2}. However this
calculation would
yield an expression that looks very different from that of \fullref{thm:monx2}, much like the two monodromy factorization given in
\fullref{prop:fuller} look very different, and comparing the
two canonical pencils would be very difficult. On the other hand,
if we can place the Lagrangian torus $T\subset X_2$ in standard position
with respect to the canonical pencil studied in \fullref{s:monx2}, then
\fullref{thm:horilu} and \fullref{prop:partialconj} allow us
to determine immediately the monodromy of a particular symplectic Lefschetz
pencil on $X_1$. We will then prove in \fullref{s:jhol} that this symplectic
Lefschetz pencil is isotopic to a generic pencil of holomorphic curves
in the canonical linear system. For now, our main result is the following:

\begin{theorem}\label{thm:monx1}
$X_1$ carries a symplectic Lefschetz pencil $\tilde{f}'$
whose fibers represent the
canonical class and whose monodromy is described by the factorization
\begin{multline*}
\prod_{i=1}^2\Biggl[\,\prod_{j=1}^5\zeta'_{i,j}\cdot\prod_{j=1}^5\xi'_{i,j}\cdot
\delta_i^+\cdot\delta_i^-\cdot(\bar\tau_1\cdot\bar\tau_2\cdot\bar\tau_3\cdot\bar\tau')^5
\Biggr]
\cdot\\
\cdot
\prod_{i=3}^6\Biggl[\,\prod_{j=1}^5\zeta_{i,j}\cdot\prod_{j=1}^5\xi_{i,j}\cdot
\delta_i^+\cdot\delta_i^-\cdot(\bar\tau_1\cdot\bar\tau_2\cdot\bar\tau_3\cdot\bar\tau_4)^5
\Biggr]
\cdot \delta_7^+\cdot\delta_7^-\cdot\delta_8^+\cdot\delta_8^-
\end{multline*}
in $\Map_{17,16}$,
where $\zeta_{i,j}$, $\xi_{i,j}$, $\delta_i^\pm$, and $\bar\tau_i$ are the
Dehn twists along the preimages of the arcs and curves
represented in \fullref{fig:monx2}, and
$\zeta'_{i,j}=\phi\zeta_{i,j}\phi^{-1}$,
$\xi'_{i,j}=\phi\xi_{i,j}\phi^{-1}$, $\bar\tau'=\phi\bar\tau_4\phi^{-1}$,
where $\phi$ is the Dehn twist along the preimage of the line segment 
joining the two leftmost branch points $\alpha$ and $\tilde\alpha_5$.
\end{theorem}

\proof We consider the loop $\gamma\subset\CP^1$ and the annulus $A$
introduced in \fullref{ss:horilu}.
Since by construction $A$ lies away from the zero section of
$\FF_6$, in a neighborhood of $A$ the pencil 
$f\co \FF_6\setminus \{p_1,\dots,p_8\}\to\CP^1$ introduced in Section 3.1 is very
close to the standard projection $pr\co \FF_6\to\CP^1$. In particular,
even though $f$ does not map $A$ to $\gamma$, there is
a nearby annulus $A'\subset\FF_6$ with boundary on $\Delta_\infty\cup C_2$
and with the property that $f(A')=\gamma$. Up to a small exact perturbation of
the symplectic form we can assume that $A'$ is Lagrangian.

Recall that by construction the loop $\gamma$
bounds a disc $\DD_-$ containing the points $q_1,z_1,r_1,q_2,z_2,r_2$.
Hence, the monodromy of $f$
along $\gamma$ (as an element of $\Map_{0,8}$) is the product of two
boundary twists (at the critical values $z_1$ and $z_2$), 
while the braid monodromy of $\Delta_\infty\cup
C_2$ along $\gamma$ is the product of the contributions from the points
inside $\DD_-$; by the calculations in \fullref{ss:calcx2}, this is
the braid which moves the points
$\alpha,\tilde{\alpha}_5,\dots,\tilde{\alpha}_1$ counterclockwise around
the points $\tilde{q}_{1,1},\dots,\tilde{q}_{1,5},z_1,\tilde{q}_{2,1},\dots,
\tilde{q}_{2,5},z_2$ by 360 degrees.
The annulus $A'$ intersects the fiber of $f$ above any point of $\gamma$ in
an arc $\eta'$ which is isotopic to a straight line segment joining the
two points of $\Delta_\infty\cup C_2$ labelled $\alpha$ and
$\tilde{\alpha}_5$; as expected, the monodromy along $\gamma$ preserves
the arc $\eta'$.

Braiding $\Delta_\infty\cup C_2$ along the annulus $A'$ yields a symplectic
curve in $\FF_6$ which is a small isotopic perturbation of the curve $C'$
considered in \fullref{ss:horilu}; in fact, for all practical purposes
we can assume that this is the same curve (for example \fullref{prop:isotopy} clearly still holds), and so we again denote it by
$C'$. By the
argument in \fullref{ss:luttinger}, the double cover $X'$ of $\FF_6$
branched along the curve $C'$ comes equipped with a
symplectic Lefschetz pencil $\tilde{f}'$ with fibers of genus 17, obtained
as a partial twisting (along the torus $T'=\pi_2^{-1}(A')$) of the pencil
$\tilde{f}_2=f\circ \pi_2$ described in \fullref{s:monx2}.

By \fullref{prop:partialconj}, the monodromy of $\tilde{f}'$ is
obtained from that of $\tilde{f}_2$ by conjugating the monodromy around each
critical value inside the disc $\DD_-$ by the Dehn twist along the
preimage of the arc $\eta'$. (Strictly speaking, we have to assume that
the loop $\gamma\subset \CP^1$ has been chosen in such a way that the base
point used for monodromy calculations in \fullref{s:monx2} lies in
$\DD_+$ and close to its boundary; however, it is easy to ensure that
this is the case). The factors that need to be conjugated in the expression
of \fullref{thm:monx2} are those corresponding to the monodromy near
the points $q_1,z_1,r_1,q_2,r_2,z_2$, namely the first 64 factors (those
corresponding to $i=1$ or $i=2$ in the product), and the conjugating 
Dehn twist is precisely $\phi$. Hence, we obtain the expression given in
the statement of \fullref{thm:monx1}.

To complete the argument, we only need to show that the
pencil $\tilde{f}'$ on $X'$ can be viewed as a symplectic
Lefschetz pencil on $X_1$ whose fibers represent the canonical class.
Indeed, by \fullref{prop:isotopy} there exists a symplectomorphism
$\psi\co \CP^1\times\CP^1\to \FF_6$ such that $\psi(C_1)=C'$.
Composing the pencil $f\co \FF_6\setminus \{p_1,\dots,p_8\}\to\CP^1$ with
$\psi$, we obtain a symplectic Lefschetz pencil on
$\CP^1\times \CP^1$, whose fibers are symplectic curves representing
the homology class $(1,4)$. Moreover, the fibers of $f\circ\psi$ intersect 
$C_1$ transversely and positively
except at isolated nondegenerate tangency points. By construction, the
symplectic Lefschetz pencil $f\circ\psi\circ \pi_1$ on $X_1$ is isomorphic
to the symplectic Lefschetz pencil $\tilde{f}'$ on $X'$ (the
isomorphism is given by the symplectomorphism from $X_1$ to $X'$
obtained by lifting $\psi$ to the double covers); and its fibers represent
the homology class $\pi_1^*(1,4)=K_{X_1}$.
\endproof

\section{Pencils of pseudo-holomorphic spheres in $\CP^1\times\CP^1$}\label{s:jhol}

Our goal in this section is to compare the symplectic Lefschetz pencil
$\tilde{f}'$ described in \fullref{thm:monx1} with a generic pencil
$\tilde{f}_1$ of holomorphic curves in the linear system $|K_{X_1}|$.
We claim:

\begin{theorem}\label{thm:isopen}
The Lefschetz pencils $\tilde{f}'$ and $\tilde{f}_1$ are isomorphic.
\end{theorem}

Recall from \fullref{s:monx1} that the Lefschetz pencil $\tilde{f}'$ is
constructed as follows. Consider the curve $C'\subset\FF_6$ obtained
by twisting $\Delta_\infty\cup C_2$ along the Lagrangian annulus $A'$,
and a pencil $f\co \FF_6\setminus \{p_1,\dots,p_8\}\to\CP^1$ of curves in
the linear system $|\Delta_0+F|$. The symplectic curve $C'$ intersects the
fibers of $f$ positively and transversely except at isolated nondegenerate
tangency points (which all lie away from $A'$), and the pencil $\tilde{f}'$
is obtained by lifting $f$ via the double cover $\pi'\co X_1\simeq X'\to \FF_6$ branched
along $C'$. Since $\CP^1\times\CP^1$ and $\FF_6$
(with the chosen K\"ahler forms in the classes $(1,4)$ and $[\Delta_0]+[F]$)
are symplectomorphic, we can also view $C'$ as a symplectic curve in $\CP^1\times
\CP^1$ representing the homology class $(6,12)$, and $f$ as a symplectic
Lefschetz pencil on $\CP^1\times\CP^1$ whose fibers represent the homology
class $(1,4)$.

The pencil $\tilde{f}_1$ can also be constructed in a similar manner, by
considering a pencil of algebraic curves in the class $(1,4)$ on
$\CP^1\times\CP^1$ whose fibers intersect the algebraic curve $C_1$
transversely except at isolated nondegenerate tangency points, and lifting
it via the double cover $\pi_1\co X_1\to\CP^1\times\CP^1$ branched along $C_1$.
Moreover, recall that Siebert and Tian's isotopy result (\fullref{thm:ST}) shows the existence of a symplectic isotopy between the
curves $C'$ and $C_1$, ie, a
continuous one-parameter family of symplectic curves $C_t\subset
\CP^1\times\CP^1$, $t\in [0,1]$, such that $C_t$ equals $C'$ for $t=0$ and 
$C_1$ for $t=1$. With this understood, \fullref{thm:isopen} is an
immediate corollary of the following statement:

\begin{proposition}\label{prop:isopen}
There exists a continuous family of symplectic Lefschetz pencils $f_t$,
$t\in [0,1]$ on $\CP^1\times\CP^1$ such that $f_0=f$, $f_1$ is a pencil of
algebraic curves, and for all $t$ the curve $C_t$ intersects the fibers of
$f_t$ positively and transversely except at isolated nondegenerate
tangency points which lie in distinct smooth fibers of $f_t$.
\end{proposition}

In fact, we will equip $\CP^1\times\CP^1$ with a family of
almost-complex structures $J_t$, $t\in [0,1]$, tamed by the fixed
symplectic form $\omega$, and work with pseudoholomorphic curves.
We start with:

\begin{lemma}
$\FF_6\simeq \CP^1\times\CP^1$ carries an almost-complex structure
$J_0$ tamed by $\omega$ and such that the curve $C'$ and the fibers of
the pencil $f$ are $J_0$--holomorphic.
\end{lemma}

\proof
By construction the fibers of $f$ are holomorphic with respect to the
standard complex structure $\JJ$ on $\FF_6$, and so is the curve
$C'$ outside of a neighborhood of the annulus $A'$. Hence we only need to
modify $\JJ$ in a neighborhood of $A'$ in order to make $C'$
pseudoholomorphic.

Over a neighborhood $U$ of $A'$ (in which $C'$ is transverse to the fibers
of $f$, and outside of which $C'$ coincides with the holomorphic curve 
$\Delta_\infty\cup C_2$), we can decompose the tangent bundle to $\FF_6$
into a direct sum $T_1\oplus T_2$, where $T_1$ is the tangent space to
the fiber of $f$ and $T_2$ is its symplectic orthogonal. This splitting
is preserved by $\JJ$, and choosing orthonormal bases of $T_1$ and $T_2$
for the metric induced by $\JJ$ and $\omega$, we have
$$\JJ=\begin{pmatrix}j_0&0\\0&j_0\end{pmatrix},\quad\mbox{where}\ 
j_0=\begin{pmatrix}0&\!-1\\1&\,\,\,0\end{pmatrix}.$$
At any point $p\in C'\cap U$, the transversality of $C'$ to the fibers of $f$
implies that we can view the tangent space $T_pC'$ as the graph of a linear
map $h\co T_2\to T_1$; the fact that $C'$ is a symplectic curve means that
$\det(h)>-1$ (with respect to the area forms induced by $\omega$ on
$T_1$ and $T_2$). The almost-complex structure
$$J=\begin{pmatrix}1&h\\0&1\end{pmatrix}\begin{pmatrix}j_0&0\\0&j_0
\end{pmatrix}\begin{pmatrix}1&h\\0&1\end{pmatrix}^{-1}=\begin{pmatrix}
j_0&hj_0-j_0h \\ 0&j_0\end{pmatrix}$$
preserves $T_1$ and $T_pC'$. Given a vector $(X,Y)\in T_1\oplus T_2$, we have
\begin{align*}\omega((X,Y),J(X,Y))&=
\omega\bigl((X,Y),(j_0X+(hj_0-j_0h)Y,j_0Y)\bigr)\\
&=\langle X,X\rangle-\langle X,(j_0hj_0+h)Y\rangle+\langle Y,Y\rangle,
\end{align*}
where $\langle\cdot,\cdot\rangle$ is the metric induced by $\JJ$ and
$\omega$. However, decomposing $h$ into its complex linear and antilinear
parts $h^{1,0}$ and $h^{0,1}$ (with respect to $j_0$), we have $\det(h)=
\|h^{1,0}\|^2-\|h^{0,1}\|^2$, so the norm of $(j_0hj_0+h)Y=2h^{0,1}Y$ is
less than twice the norm of $Y$. This implies that $J$ is tamed by $\omega$.

We have therefore obtained an almost-complex structure with the desired
properties at every point of $C'\cap U$; moreover, near the boundary of $U$ the
curve $C'$ is $\JJ$--holomorphic and hence $h$ is complex linear, so that
$J$ coincides with $\JJ$. We can extend this construction to a tubular
neighborhood of $C'\cap U$ by choosing a suitable extension of $h$
(preserving the condition $\det(h)>-1$ and the complex linearity near the
boundary of $U$); patching this together with $\JJ$ by means of a suitable
cut-off function, we obtain a
globally defined almost-complex structure with the desired properties.
\endproof

Starting from the $J_0$--holomorphic curve $C'$, the method used by
Siebert and Tian to prove symplectic isotopy \cite{ST} yields a family of
$\omega$--tame almost-complex structures $J_t$, $t\in [0,1]$ on
$\CP^1\times \CP^1$, with $J_1$ equal to the standard (product) complex
structure, and a family of smooth $J_t$--holomorphic curves $C_t$ realizing
the isotopy between $C'$ and $C_1$.

At this point, we need to review some standard results about pseudoholomorphic
spheres in $\CP^1\times\CP^1$. As observed by Hofer--Lizan--Sikorav~\cite{HLS},
the linearized $\bar\partial$--operator is always surjective for embedded
pseudoholomorphic spheres of self-intersection number at least $-1$ in
an almost-complex 4--manifold (see also \cite[Lemma 3.3.3]{McS}); this
property is sometimes called {\it automatic regularity}. Hence, the moduli
spaces of $J_t$--holomorphic spheres (ie, embedded irreducible 
$J_t$--holomorphic curves of genus 0) considered below are always smooth
manifolds of the expected dimensions (provided they are non-empty).
This implies:

\begin{lemma} \label{l:jspheres}
Let $J$ be any $\omega$--tame almost-complex structure on $\CP^1\times\CP^1$.

{\rm(i)}\qua Any point $p\in\CP^1\times\CP^1$ lies on a unique 
$J$--holomorphic sphere representing the homology class $(0,1)$, which we
call the {\em $J$--fiber} through~$p$.

{\rm(ii)}\qua Given an integer $k\ge 0$ and $2k+1$ distinct points
$p_1,\dots,p_{2k+1}$ in $\CP^1\times\CP^1$, there exists at most one
$J$--holomorphic sphere representing the homology class $(1,k)$ and passing
through the points $p_1,\dots,p_{2k+1}$.

{\rm(iii)}\qua Fix $k\ge 0$, and assume that $\CP^1\times\CP^1$ contains no
$J$--holomorphic spheres in the homology class $(1,j)$ for any $j<-k$.
Let $p_1,\dots,p_{2k+1}$ be distinct points in $\CP^1\times
\CP^1$, such that no two of them lie in the same $J$--fiber.
Assume moreover that for all $-k\le j<k$, no $j+k+1$ of the points $p_1,\dots,p_{2k+1}$
lie on the same $J$--holomorphic sphere in the homology class $(1,j)$.
Then there exists a unique $J$--holomorphic sphere representing the
homology class $(1,k)$ and passing through the points $p_1,\dots,p_{2k+1}$.
\end{lemma}

\proof The first statement is classical and due to Gromov
\cite[Theorem 0.2.A]{Gr}. The second statement is also classical and
follows from positivity of intersections: if two $J$--holomorphic
curves in the class $(1,k)$ intersect in $2k+1$ distinct points then
they must share a component; since we assume irreducibility, they must
be equal.

To prove the third statement, we use the fact that the
Gromov--Witten invariant which counts pseudoholomorphic curves of genus $0$
in the class $(1,k)$ passing through $2k+1$ points is non-zero. Indeed,
when $J$ is the standard
complex structure and $p_1,\dots,p_{2k+1}$ are generic, the $2k+1$ incidence
conditions determine a one-dimen\-sional linear subspace in the vector space
$H^0(\OO_{\PP^1\times \PP^1}(1,k))$, ie, there is a unique algebraic curve
through the given points. For a generic choice of the points this curve is
smooth and hence automatically regular; since the complex structure is
integrable, its contribution to the Gromov--Witten invariant is~$1$.

Returning to the case of arbitrary $J$, this implies
the existence of a (possibly singular) $J$--holomorphic curve of genus 0
through $p_1,\dots,p_{2k+1}$ in the homology class $(1,k)$. We claim
that the assumptions on $p_1,\dots,p_{2k+1}$ imply smoothness.
Indeed, if the curve is not smooth then it must be reducible and a union of
smoothly embedded $J$--holomorphic spheres (this follows
eg from the adjunction formula). However, by positivity of intersection
with the $J$--fibers, every irreducible $J$--holomorphic curve must represent a homology
class of the form $(a,b)$ with $a\ge 0$; and if $a=0$ then necessarily
$b=1$ (positivity of area implies $b\ge 1$, and adjunction implies $b=1$).
Therefore, our curve must be the union of a $J$--holomorphic sphere
representing the homology class $(1,j)$ for some
integer $-k\le j<k$, and $k-j$ fibers. However, by assumption each of the
$k-j$ fibers contains at most one of the points $p_1,\dots,p_{2k+1}$, and
the component representing the class $(1,j)$ passes through at most $j+k$
of them. This yields a contradiction.
\endproof

\begin{definition}
We say that a configuration of 8 distinct points $p_1,\dots,p_8\in
\CP^1\times\CP^1$ is {\em $J$--regular} if the following conditions hold:
none of them lie on a $J$--holomorphic sphere in the homology class $(1,j)$
for some $j<0$; no two of them lie in a same $J$--fiber; and for all $j\in
\{0,1,2,3\}$, no $2j+2$ of them lie on a same $J$--holomorphic sphere in the
homology class $(1,j)$.
\end{definition}

\begin{lemma} \label{l:jreg}
For any $\omega$--tame almost-complex structure, the set of $J$--regular
configurations is a connected open subset of $(\CP^1\times\CP^1)^8$ whose
complement has real codimension 2.
\end{lemma}

\proof 
By positivity of intersections, there is at most one $J$--holomorphic sphere
representing a homology class of the form $(1,j)$ with $j<0$; configurations
containing a point on such a sphere therefore form a codimension 2 subset.
The moduli space of $J$--fibers has real dimension 2, so the space of
configurations of two points on a same $J$--fiber has real dimension 6, and
the space of configurations of 8 points of which two lie on a same $J$--fiber
has real dimension 30, ie, codimension 2. Similarly, for $0\le j\le 3$,
automatic regularity implies that the moduli space of $J$--holomorphic
spheres in the class $(1,j)$ has real dimension $4j+2$, the space of
configurations of $2j+2$ points on such a sphere has real dimension $8j+6$,
and the space of configurations of 8 points of which $2j+2$ lie on such a
sphere again has codimension 2.
\endproof

\begin{proposition} \label{prop:jpen}
Let $(p_1,\dots,p_8)$ be a $J$--regular configuration of points in
$\CP^1\times\CP^1$. Then the family of all $J$--holomorphic curves of
genus $0$ which represent the homology class $(1,4)$ and pass through
$p_1,\dots,p_8$ forms a Lefschetz pencil. Moreover, the 8 singular fibers
of this pencil are reducible $J$--holomorphic curves consisting
of the $J$--fiber through some $p_i$ $(1\le i\le 8)$ and the $J$--holomorphic
sphere in the class $(1,3)$ through all $p_j$, $j\neq i$.
\end{proposition}

\proof Since $[\omega]$ is Poincar\'e dual to $(1,4)$, there are no
symplectic spheres in the homology classes $(1,j)$ for $j\le -4$.
Hence, \fullref{l:jspheres} {\rm(iii)} implies the existence of a unique
$J$--holomorphic sphere in the homology class $(1,3)$ through any 7 of the
points $p_1,\dots,p_8$.

Consider any point $p\in\CP^1\times\CP^1\setminus\{p_1,\dots,p_8\}$.
There are two cases.
If $p$ does not lie on any of the
$J$--fibers through $p_1,\dots,p_8$, nor on any of the $J$--holomorphic
spheres in the homology class $(1,3)$ through seven of these points,
then one easily checks that the 9 points $p_1,\dots,p_8,p$ satisfy the
assumptions of \fullref{l:jspheres} {\rm(iii)}, and
so they lie on a unique
$J$--holomorphic sphere in the homology class $(1,4)$.
Otherwise, $p_1,\dots,p_8,p$ lie on a reducible $J$--holomorphic curve of
the type described in the statement of the proposition, and by positivity
of intersections no other $J$--holomorphic curve
representing the homology class $(1,4)$ can pass through these 9 points.

To see that this family of $J$--holomorphic curves is parameterized by
$\CP^1$, consider a $J$--fiber $F$ (not passing through any $p_i$), and
observe that each curve in the family intersects $F$ transversely in a
single point, and conversely through any point of $F$ there is a single
curve in the family. Hence, we can define a map from
$\CP^1\times\CP^1\setminus \{p_1,\dots,p_8\}$ to $\CP^1\simeq F$ by
mapping each point $p$ to the point where the curve through
$p_1,\dots,p_8,p$ intersects $F$.
The fact that this is a Lefschetz pencil follows
from automatic regularity and standard arguments about deformations
of $J$--holomorphic curves; in particular, the structure of the moduli
space near the nodal curves
follows from gluing arguments (see eg \cite[Corollary 2]{Sik}).
In fact, the key ingredient is again automatic regularity, which implies
that the local behavior of families of $J$--holomorphic spheres is the same
as in the usual holomorphic case.
\endproof   

\proof[Proof of \fullref{prop:isopen}]
Consider the family of almost-complex structures $J_t$ and pseudoholomorphic
curves $C_t$ introduced above.
By \fullref{l:jreg}, we can find a continuous family of $J_t$--regular
configurations $(p_{1,t},\dots,p_{8,t})$, $t\in [0,1]$, starting from the
configuration chosen in \fullref{ss:x2notation} at $t=0$ (which is easily
seen to be $J_0$--regular, recalling that $\Delta_0\subset\FF_6$ represents
the class $(1,3)$ in $\CP^1\times\CP^1$). Moreover, we can choose these
points in such a way that, for all $1\le i\le 8$, the point $p_{i,t}$ does
not lie on $C_t$, and the point where the $J_t$--fiber through $p_{i,t}$
intersects the $J_t$--holomorphic sphere in the class $(1,3)$ through
all $p_{j,t}$, $j\neq i$, does not lie on $C_t$ either.

By \fullref{prop:jpen}, for each $t$ the points
$(p_{1,t},\dots,p_{8,t})$ determine a pencil $f_t$ of
$J_t$--holomorphic curves in the homology class $(1,4)$.  Since the
curve $C_t$ is irreducible, every intersection of $C_t$ with a fiber
of $f_t$ has a finite positive multiplicity; see eg \cite[Appendix
E]{McS} for a detailed discussion of this fact, which follows from
Micallef and White's result about the local structure of
pseudoholomorphic curves \cite{MiWh}. Moreover, the restriction of
$f_t$ to $C_t$ is an open mapping $\phi_t\co C_t\to\CP^1$ of degree 36,
whose critical points are precisely the non-transverse intersections
between $C_t$ and fibers of $f_t$.

The preimage by $\phi_t$ of a small disc in $\CP^1$ centered
at a critical value $z_{cr}$ of $\phi_t$ (chosen generically so its
boundary is transverse to $\phi_t$) consists of at most
$|\phi_t^{-1}(z_{cr})|\le 35$ components, each of which has Euler
characteristic at most 1. Hence, arguing as in the classical Hurwitz
formula for branched covers, we conclude
that $\phi_t$ has at most 180 critical points; in particular, the points
where $C_t$ is tangent to the fibers of $f_t$ are isolated.

We conclude that $\phi_t\co C_t\to\CP^1$ is a topological branched covering.
After modifying $\phi_t$ by a $C^1$--small perturbation supported in a
neighborhood of its critical points, we can assume that the critical points
of $\phi_t$ are all non-degenerate, and that the corresponding
critical values are distinct from each other and from the critical values
of $f_t$. Using suitable cut-off functions, this modification of $\phi_t$
can be extended to a $C^1$--small perturbation of $f_t$ supported in a
neighborhood of the critical points of $\phi_t$, preserving the property
that the fibers of the pencil intersect $C_t$ positively.
The fibers of the perturbed pencil are no
longer $J_t$--holomorphic, but they can still be assumed to be symplectic.
It is moreover clear that this perturbation argument can be carried out
in a manner such that the perturbations depend continuously on $t\in [0,1]$.
\endproof

\section{Comparing the canonical pencils}\label{s:mainthm}

We now have all the necessary ingredients to compare generic pencils of
curves in the canonical linear systems on $X_1$ and $X_2$. In particular,
\fullref{thm:main1} follows directly from \fullref{thm:monx2},
\fullref{thm:monx1}, and \fullref{thm:isopen}. Moreover, in order
to compare the monodromy groups and prove \fullref{thm:subgp}, it is
enough to prove that the conjugating element $\phi$ belongs to the monodromy
subgroups of both pencils. Namely, we have to prove:

\begin{proposition}\label{prop:subgp}
With the notations of \fullref{thm:main1}, the Dehn twist $\phi$ belongs
to the subgroup $G_2$ of \,$\Map_{17,16}$ generated by $t_1,\dots,t_{196}$,
and it also belongs to the subgroup $G_1$ generated by $\phi t_1\phi^{-1},
\dots,\phi t_{64}\phi^{-1},t_{65},\dots,t_{196}$. Therefore, $G_1=G_2$.
\end{proposition}

It is easy to prove that $\phi\in G_1$. Indeed, recall that Dehn
twists along simple closed curves that intersect transversely
once satisfy the relation $t_at_bt_a=t_bt_at_b$. Therefore, with the
notations of \fullref{thm:monx1}, we have
 $\bar{\tau}'=\phi\bar{\tau}_4\phi^{-1}=\bar{\tau}_4^{-1}
\phi\bar{\tau}_4$, so $\phi=\bar{\tau}_4\bar{\tau}'\bar{\tau}_4^{-1}$
obviously belongs to the subgroup $G_1$ generated by the various Dehn twists
appearing in \fullref{thm:monx1}.

The argument for $X_2$ is more subtle, since a quick inspection of the
factors in \fullref{thm:monx2} does not suggest any obvious reason
why $\phi$ should belong to the monodromy group. We use
the same notations as in \fullref{ss:monx2}; in particular, we consider
Dehn twists in $\Map_{17,16}$ which are obtained by lifting
half-twists or Dehn twists via the double cover 
$\pi_2\co \tilde\Sigma\to\Sigma$, where $\Sigma=\CP^1\setminus\{z_1,\dots,z_8\}$
and the 36 branch points of $\pi_2$ are labelled as in \fullref{ss:monx2}.

\begin{figure}[ht!]
\centering
\includegraphics{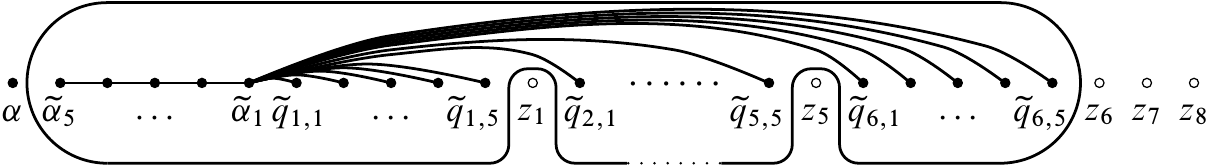}
\caption{The disc $\DD_0$ and the tree $T_0$}\label{fig:DD0}
\end{figure}

\begin{lemma} \label{l:b35}
The subgroup of\/ $\Map_{17,16}$ generated by the Dehn twists $\zeta_{i,j}$
and $\bar\tau_j$ is the image by the lifting homomorphism of the braid group
$B_{35}$ consisting of all braids supported in a disc $\DD_0\subset\Sigma$
which contains the $35$ branch points
$\tilde\alpha_j$ and $\tilde{q}_{i,j}$ as well as arcs connecting them within
the upper half-plane, as shown on \fullref{fig:DD0}.
\end{lemma}

\proof Conjugating $\zeta_{i,j}$ by $\bar\tau_{j-1}\dots\bar\tau_1$,
we obtain the lift of a half-twist that exchanges the branch points
$\tilde\alpha_1$ and $\tilde{q}_{i,j}$ along an arc contained in the upper
half-plane. These 30 arcs together with the 4 arcs supporting the half-twists
which lift to $\bar\tau_1,\dots,\bar\tau_4$ form an embedded tree 
$T_0\subset\Sigma$ (see \fullref{fig:DD0}). It is well-known that these
half-twists generate the braid group $B_{35}$ (for example, further
conjugations yield half-twists whose supporting arcs form a linear chain
as in Artin's standard set of generators). \endproof

\begin{definition}
The {\em upper envelope} of a subset $S$ of 
$S_0=\{\alpha,\,\tilde{\alpha}_j,\,\tilde{q}_{i,j},\ 
i=1\dots 6,$ $j=1\dots 5\}\cup\{z_i,\ 1\le i\le 8\}$
is the simple closed curve $c(S)\subset\Sigma$ which
bounds a disc containing the points of $S$ as well as arcs connecting
them within the upper half-plane, but not any points of $S_0\setminus S$.
We denote by $\delta(S)^\pm$ the two lifts of the Dehn twist along $c(S)$.
\end{definition}

For example, the factors $\delta_i^{\pm}$ in \fullref{thm:monx2} are 
the lifts of the Dehn twists along the upper envelopes
of the sets $\{\alpha,\tilde{\alpha}_5,\dots,\tilde{\alpha}_1,z_i\}$.
With this terminology, the following result is an immediate corollary of
\fullref{l:b35}:

\begin{lemma}\label{l:permut}
If the monodromy group $G_2$
contains $\delta(S)^\pm$ for some set $S$, then it also contains
$\delta(S')^\pm$ for any set $S'$ which is the image of $S$ by a
permutation of $S_0$ fixing the elements $\alpha,z_1,\dots,z_8$. 
\end{lemma}

\begin{figure}[ht]
\centering
\includegraphics{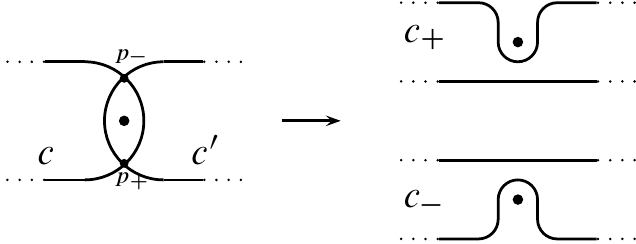}
\caption{Summing simple closed curves}\label{fig:sum}
\end{figure}

In particular, $G_2$ contains $\delta(S)^\pm$ for any 7--element set $S$
which contains $\alpha$ and exactly one $z_i$. Our next observation is
the following:

\begin{lemma}\label{l:sum}
Let $c,c'$ be two simple closed curves in $\Sigma$, intersecting in two points
as in \fullref{fig:sum}, left, so one of the regions delimited by $c\cup c'$
contains a single branch point of $\pi_2$. Assume that $G_2$ contains the Dehn twists
along both lifts of $c$ and $c'$. Then $G_2$ also contains the
Dehn twists along both lifts of the loops $c_\pm$
obtained by ``summing'' $c$ and $c'$ (\fullref{fig:sum}, right).
\end{lemma}

\proof Let $\tilde{c}$ and $\tilde{c}'$ be arbitrary lifts of $c$ and $c'$.
The loops $\tilde{c}$ and $\tilde{c}'$ intersect
only once, at a point of $\pi_2^{-1}(c\cap c')$ which depends on the chosen
lifts. Hence, denoting by $\tilde{t}$ and $\tilde{t}'$ the corresponding
Dehn twists, the conjugate of $\tilde{t}'$ by $\tilde{t}$ is the
Dehn twist along the connected sum of $\tilde{c}$ with $\tilde{c}'$.
Assume that the intersection of $\tilde{c}$ with $\tilde{c}'$ lies above
$p_+$: then $\tilde{t}\tilde{t}'\tilde{t}^{-1}$ 
is the Dehn twist along $\tilde{t}^{-1}(c')$, which is a
lift of $c_+$ since the lift of a loop that encircles the
branch point twice is contractible.
Similarly, if the intersection of $\tilde{c}$ and $\tilde{c}'$
lies above $p_-$, then $\tilde{t}^{-1}\tilde{t}'\tilde{t}$ is the Dehn twist
along a lift of $c_-$. Considering the
four possible choices for $(\tilde{c},\tilde{c}')$, we obtain both lifts
of $c_+$ and both lifts of $c_-$. \endproof

In particular, let $S$ and $S'$ be two subsets of $S_0$ such that 
$S\cap S'$ consists of a single element $s$ which is not one of the $z_i$.
Assume moreover that, with respect to the natural ordering of the elements
of $S_0$ induced by their positions along the real axis, the following 
conditions are satisfied:
\begin{itemize}
\item $s\in \{\inf(S'),\sup(S')\}$;
\item $S\cap [\inf(S'),\sup(S')]=\{s\}$.
\end{itemize}
In this situation, if the monodromy group $G_2$ contains $\delta(S)^\pm$ and
$\delta(S')^\pm$, then by \fullref{l:sum} it also contains
$\delta(\bar{S})^\pm$ where $\bar{S}=(S\cup S')\setminus \{s\}$.

Applying this argument repeatedly to specific subsets of $S_0$ satisfying
these conditions, and combining with \fullref{l:permut},
we conclude that $G_2$ contains
$\delta(S)^{\pm}$ whenever $S$ satisfies one of the following conditions:

\begin{itemize}
\item $S$ contains 7 elements, among which one $z_i$, and $\alpha\in S$;
\item $S$ contains 12 elements, among which two $z_i$, and $\alpha\not\in S$;
\item $S$ contains 17 elements, among which three $z_i$, and $\alpha\in S$;
\item $S$ contains 22 elements, among which four $z_i$, and $\alpha\not\in S$;
\item $S$ contains 27 elements, among which five $z_i$, and $\alpha\in S$;
\item $S$ contains 32 elements, among which six $z_i$, and $\alpha\not\in S$;
\item $S$ contains 37 elements, among which seven $z_i$, and $\alpha\in S$;
\item $S$ contains 42 elements, among which eight $z_i$, and $\alpha\not\in S$.
\end{itemize}

\noindent
In particular, taking $S_\phi=S_0\setminus \{\alpha,\tilde\alpha_5\}$ (which
contains 42 elements, including all the $z_i$) we obtain that 
$\delta(S_\phi)^\pm\in G_2$; however it is easy to see that 
$\delta(S_\phi)^\pm=\phi$, which completes the proof of \fullref{prop:subgp}.

\section{Pluricanonical pencils and degree doubling}\label{ss:doubling}

As mentioned in \fullref{ss:slf}, comparing Horikawa surfaces requires
understanding not only pencils of curves in the canonical linear system,
but also in the pluricanonical linear systems $|kK_{X_i}|$, $k\gg 1$.
This can be achieved by means of a {\it degree doubling} procedure, which
describes the topology of a pencil of curves in the linear system $|2kK_{X_i}|$
in terms of that of a pencil in the linear system $|kK_{X_i}|$. This idea, which
goes back to Donaldson \cite{Do1}, has been explored in greater detail in \cite{Sm}
and \cite{AK}. We start by giving an outline of the relevant material
in those two papers.

Consider a Lefschetz pencil $f=(s_0/s_1)$ defined by two generic sections $s_0, s_1$ of a
sufficiently positive line bundle $L^{\otimes k}$ (either holomorphic, or
approximately holomorphic in the sense of Donaldson \cite{Do2}). As observed
in \cite{Sm}, the two sections
$s_0^2$ and $s_0s_1$ of $L^{\otimes 2k}$ define a (highly
non-generic) pencil of reducible nodal curves $\{s_0^2-\alpha s_0s_1=0
\}_{\alpha\in\CP^1}$, obtained by adding in the zero set of $s_0$ to each of
the curves in the original pencil. A generic pencil can be obtained by 
choosing small perturbations $\epsilon_0,\epsilon_1$ 
(sections of $L^{\otimes 2k}$) and considering the sections
$s_0^2+\epsilon_0$ and $s_0s_1+\epsilon_1$ instead. It is then easy to see
that the generic fiber of this pencil is obtained by forming the connected
sum of two generic fibers of $f$ (smoothing the intersections at the base
points), and that the critical points and vanishing cycles which occur away
from the zero set of $s_0$ are in one-to-one correspondence with those of $f$
\cite{Sm}. If $f$ has fiber genus $g$ and $n$ base points, then the
doubled pencil has fiber genus $\bar{g}=2g+n-1$ and $\bar{n}=4n$ base points.
Its monodromy consists of: (1) the
image of the monodromy of $f$ under a natural embedding $\Map_{g,n}\hookrightarrow
\Map_{\bar{g},\bar{n}}$ induced by viewing the fiber of $f$ (minus a
neighborhood of its base points) as a subset of the new fiber, and (2)
contributions from a neighborhood of the zero set of $s_0$ \cite{Sm}.

The work in \cite{AK} aims to turn these considerations into an explicit
formula for the monodromy of the doubled pencil. The starting point is that
complex surfaces, and more generally symplectic 4--manifolds (see
\cite{AK1}), can be realized as branched covers of $\CP^2$ by choosing
three suitable sections $s_0,s_1,s_2$ of the line bundle
$L^{\otimes k}$, and considering the map $F=(s_0:s_1:s_2)\co X\to\CP^2$.
Choosing the sections generically in the holomorphic case, or in a specific
manner in the approximately holomorphic case, we can ensure that $F$ is a
branched covering with generic local models, ie, near any point of the
ramification curve $R\subset X$ it is modelled on one of the two maps
$(x,y)\mapsto (x^2,y)$ or $(x,y)\mapsto (x^3-xy,y)$; moreover the branch
curve $D=F(R)\subset\CP^2$ can be assumed to have the following properties:

\begin{enumerate}
\item the only singularities of $D$ are ordinary double points (``nodes'')
and ordinary cusps (in the approximately holomorphic case, there may be
some double points with the anti-complex orientation);
\item $F_{|R}\co R\to D$ is an immersion everywhere except at the cusps, and
one-to-one except at the nodes;
\item $(0:0:1)\not\in D$;
\item $D$ is positively transverse to the fibers of the linear projection
$\pi:$ \hbox{$(x\!:\!y\!:\!z)$} \hbox{$\mapsto (x\!:\!y)$} from $\CP^2$ to $\CP^1$, except at isolated
nondegenerate tangency points (distinct from the cusps and nodes);
\item the cusps, nodes, and tangency points lie in different fibers of $\pi$.
\end{enumerate}

\noindent
Then the composition $f=\pi\circ F=(s_0/s_1)$ is a Lefschetz pencil.
The singular fibers of $f$ are the preimages by $F$ of those fibers of $\pi$
which are tangent to $D$, its base points are the preimages of
$(0\!:\!0\!:\!1)$, and its monodromy can be determined from the braid monodromy of $D$ using
a lifting homomorphism as in \fullref{ss:lifting}.

As observed in \cite{AK}, the composition of $F$ with a generic quadratic
map $V'_2\co \CP^2\to\CP^2$ yields a map $\bar{F}=V'_2\circ F\co X\to\CP^2$ determined
by three sections of $L^{\otimes 2k}$ which are quadratic expressions in
$s_0,s_1,s_2$; this map is again a branched covering, ramified along
$\bar{R}=R\cup F^{-1}(R')$ (where $R'$ is the ramification curve of $V'_2$). The
branch curve of $\bar F$ consists of two parts, namely $V'_2(D)$ on
one hand, and $n=\deg F$ superimposed copies of the branch curve $D'$ of
$V'_2$ on the other hand. Near a point where $F^{-1}(R')$ intersects $R$,
a local model for $\bar F$ is given by the composition of two simple
branched covers such that the branch curve of the first map is in general
position with respect to the second one, $(x,y)\mapsto (x^2+y,y)\mapsto 
(x^2+y,y^2)$. At such a point, the branch curve of $\bar{F}$ is
not immersed, and presents a self-tangency (since $V'_2(D)$ is tangent to
$D'$). So, $\bar{F}$ is not everywhere given by one of the generic
local models, and its branch curve $\bar{D}$ does not satisfy
properties (1) and (2) above: in addition to nodes and cusps, $\bar{D}$ also
has self-tangencies, and while the restriction of $\bar{F}$ to its
ramification curve is still an immersion outside the cusps and
self-tangencies, it is not generically one-to-one. Nonetheless, a small
perturbation can be added to $\bar{F}$ in order to get a covering with
generic local models, satisfing properties (1)--(5). The main idea in
\cite{AK} is that the topology of this covering (in particular the braid
monodromy of its branch curve) can be determined explicitly from that of
$F$, using the non-generic map $\bar{F}$ as an intermediate step.
 
After composing with the linear projection $\pi$, we again obtain a
Lefschetz pencil $\bar{f}$, whose monodromy can be determined by lifting the
half-twists in the braid monodromy of the branch curve. This leads to an
explicit degree doubling formula for Lefschetz pencils obtained from
sections of sufficiently positive line bundles \cite[Theorem 4, Section 4.2]{AK}:

\begin{theorem}\label{thm:k2k}
Let $f=\pi\circ F$ be a Lefschetz pencil with fiber genus $g$ and $n$ base
points, determined by two sections of $L^{\otimes k}$, where $k$ is assumed
to be sufficiently large. Let $\Phi$ be a factorization of the boundary
twist in $\Map_{g,n}$ describing the monodromy of $f$. Let $\bar{f}$ be
the Lefschetz pencil obtained by the construction described above (so
$\bar{f}$ is determined by two sections of $L^{\otimes 2k}$,
its fiber genus is $\bar{g}=2g+n-1$, and it has $\bar{n}=4n$ base points).

Then the monodromy of $\bar{f}$ can be described by the factorization
$\iota(\Phi)\cdot U_{g,n}$ in $\Map_{\bar{g},\bar{n}}$, where
$\iota\co \Map_{g,n}\hookrightarrow \Map_{\bar{g},\bar{n}}$ is a natural
embedding induced by viewing the fiber of $f$ (minus a neighborhood of its
base points) as a subset of the fiber of $\bar{f}$, and $U_{g,n}$ is an
explicitly determined
collection of $4g-4+7n$ Dehn twists in $\Map_{\bar{g},\bar{n}}$ that depends only
on $g$ and $n$ (but not on $f$).
\end{theorem}

It is not immediately clear that this approach applies to the
canonical pencils on the Horikawa surfaces $X_1$ and $X_2$. The main issue
is that the canonical pencils do not satisfy the ``large $k$'' requirement.
In particular, it is
not clear that holomorphic or approximately holomorphic perturbations
with the required properties can be constructed, and it is not clear
either that a Lefschetz pencil obtained by approximately holomorphic
techniques would be topologically equivalent to a pencil of holomorphic
curves. In fact, the linear systems $|K_{X_i}|$ and
$|2K_{X_i}|$ factor through $\CP^1\times\CP^1$ and $\FF_6$, so a generic
triple of holomorphic sections of the canonical bundle does not even
determine a map to $\CP^2$ with generic local models.

However, the features of the maps to $\CP^2$ that naturally arise in this
context are actually not an obstacle. Indeed, given a branched covering
map $F\co X\to\CP^2$, the critical points of $\pi\circ F$ are the points of
the ramification curve where the image of $dF$ is not transverse to the
fiber of $\pi$, ie, the points of $R$ which map to the vertical
tangencies of $D$; in particular, among the properties listed above, only
(3) and (4) really matter.

\begin{definition}
A branched covering map $F\co X\to\CP^2$ is {\em tame} if near any point of
the ramification curve $R\subset X$ it is modelled on one of the three
maps $(x,y)\mapsto (x^2,y)$, $(x,y)\mapsto (x^3-xy,y)$, or $(x,y)\mapsto
(x^2+y,y^2)$, and moreover the branch curve $D=F(R)\subset\CP^2$ satisfies
the following properties:
\begin{enumerate}
\item[(1$'$)] the only singularities of $D$ are ordinary double points, ordinary
cusps, and self-tangencies;
\item[(2$'$)] $F_{|R}\co R\to D$ is an immersion away from the cusps and
self-tangencies;
\item[(\,3\,\unsp)] $(0:0:1)\not\in D$;
\item[(\,4\,\unsp)] $D$ is positively transverse to the fibers of $\pi\co (x:y:z)\mapsto
(x:y)$, except at isolated nondegenerate tangency points (distinct from the
cusps, nodes, and self-tangencies).
\end{enumerate}
\end{definition}

If $F$ is a tame covering, then the composition $\pi\circ F$ is still a
Lefschetz pencil, althouh its critical points need not lie in distinct
fibers: for example, whenever a fiber of $\pi$ is tangent to a component of 
$D$ over which $F_{|R}$ is not generically one-to-one, we get a fiber
of $\pi\circ F$ with multiple nodes. With this understood, we can still
consider the individual Dehn twists obtained by
lifting the half-twists in the braid monodromy of $D$; in the case of a
component of multiplicity $\mu$, we obtain $\mu$ different Dehn twists
along disjoint simple closed curves obtained by considering appropriate
lifts of the supporting arc of the half-twist.

In the holomorphic setting, mild genericity conditions ensure that the
composition of two tame coverings is still a tame covering. In particular,
if $F\co X\to\CP^2$ is a tame covering defined by a triple of holomorphic
sections of $L^{\otimes k}$
and $V'_2\co \CP^2\to\CP^2$ is a generic quadratic map, then $V'_2\circ F$ is
still a tame covering. Moreover, considering specifically the Horikawa
surfaces $X_i$ ($i\in\{1,2\}$), the map to $\CP^2$ defined
by a generic triple of sections of the canonical bundle $K_{X_i}$
is the composition of the double covering $\pi_i$ with a generic branched
covering defined by three sections of $\mathcal{O}(1,4)$ on $\CP^1\times\CP^1$ or
$\mathcal{O}(\Delta_0+F)$ on $\FF_6$; such a map is a tame covering.

Given a tame covering map $F$, one can always modify it by an
arbitrarily small $C^\infty$ perturbation (locally
holomorphic near the vertical tangencies, cusps, and self-tangencies) in
order to obtain
a symplectic branched covering with generic local models satisfying
properties (1)--(5). More precisely, the effect of such a perturbation
is to replace each self-tangency by three cusps (replacing $(x,y)\mapsto
(x^2+y,y^2)$ by $(x,y)\mapsto (x^2+y,y^2+\epsilon x)$ for a small
nonzero $\epsilon$), and to separate each multiple component of $D$ into
distinct copies intersecting at nodes (see \cite{AK}). However, as far as
the corresponding Lefschetz pencil is concerned, the only effect of the
perturbation is to move the critical points of $\pi\circ F$ into distinct
nearby fibers; so the monodromy still consists of the same Dehn twists,
independently of the chosen perturbation.

Hence, while perturbations are needed in order to study the
effect of degree doubling on the braid monodromy of the branch curves,
which is the method used in \cite{AK} in order to derive \fullref{thm:k2k}, they are actually irrelevant as far as pencils are
concerned. In particular, in our setting the formula in \fullref{thm:k2k}
can be interpreted more directly as a relation between
the monodromies of the pencils $\pi\circ F$ and $\pi\circ V'_2\circ F$
(with the understanding that, when several critical points lie in a same
fiber, we still consider the individual Dehn twists separately).
In conclusion, we have:

\begin{proposition}
Let $F\co X\to\CP^2$ be a holomorphic map from a complex surface to $\CP^2$,
and assume that $F$ is a tame branched covering. Let $V'_2\co \CP^2\to\CP^2$
be a map defined by three generic quadratic polynomials. Then the maps
$f=\pi\circ F$ and $\bar{f}=\pi\circ V'_2\circ F$ are Lefschetz pencils,
and their monodromies
are related by the formula in \fullref{thm:k2k}.
\end{proposition}

In particular, the monodromy of a generic pencil of curves 
in the linear system $|2^m K_{X_i}|$ on the
Horikawa surface $X_i$ ($m\ge 1$, $i\in\{1,2\}$) 
consists of two ingredients:
\begin{itemize}
\item the image of the monodromy of
the canonical pencil $\tilde{f}_i$ under a natural embedding
of $\Map_{17,16}$
induced by viewing the fiber of $\tilde{f}_i$ (minus a neighborhood of
the base points) as a subset of the new fiber;
\item an
explicit collection of Dehn twists that depends on $m$ but not on $i$.
\end{itemize}
With this understood, \fullref{thm:pluri} becomes an easy corollary
of Theorems \ref{thm:main1} and \ref{thm:subgp}.

\section{Matching paths and Lagrangian spheres in $X_i$}\label{ss:matching}

Lefschetz pencils can be used to understand Lagrangian spheres in a
symplectic manifold via {\it matching paths}, an idea due to Donaldson
and Seidel (see \cite{Do1} and \cite[Section 9b]{Se}). In the
four-dimensional case, the definition is quite simple:

\begin{definition}
Let $f\co X^4\setminus B\to S^2$ be a symplectic Lefschetz pencil.
An embedded arc $\gamma\co [0,1]\to S^2$ with $\gamma^{-1}(\crit f)=\{0,1\}$
is a {\em matching path} for $f$ if the vanishing cycles associated to
the arcs $\gamma([0,\frac12])$ and $\gamma([\frac12,1])$ are isotopic
to each other inside $f^{-1}(\gamma(\frac12))\setminus B$.
\end{definition}

For example, if a same Dehn twist is repeated twice in the monodromy
factorization associated to $f$, then the ``simplest'' arc that joins the two
corresponding critical values by passing through the chosen base point is
a matching path. More
generally, matching paths arise whenever an arbitrary sequence of Hurwitz
moves leads to a factorization in which a same Dehn twist is repeated twice.

A matching path gives rise to an embedded Lagrangian sphere in $X$
(up to isotopy), obtained by joining together the two
{\it thimbles} formed by parallel transport of the vanishing cycles along
the arcs $\gamma([0,\frac12])$ and $\gamma([\frac12,1])$; see
\cite[Section 9b]{Se}. Conversely, as
observed by Donaldson, any
Lagrangian sphere can be obtained (up to isotopy) from a matching path in
a Lefschetz pencil of sufficiently high degree (see \cite{AMP} for a proof).

Matching paths can be viewed as specific elements (``figure 8 loops'')
in the kernel of the monodromy morphism $\psi\co \pi_1(\DD\setminus
\crit(f))\to \Map_{g,n}$ associated to the Lefschetz pencil, or as specific
pairs $(\gamma_+,\gamma_-)$ of conjugates of generators of
$\pi_1(\DD\setminus\crit(f))$ for which the monodromies coincide.

In the case of the canonical pencil $\tilde{f}_2$ on $X_2$, whose monodromy
has been described in \fullref{thm:monx2}, there are obvious matching
paths arising from the repeated factors $\bar\tau_i$, and slightly less
obvious matching paths arising from the fact that the conjugate of 
$\bar\tau_i$ by $\bar\tau_{i+1}$ equals the conjugate of $\bar\tau_{i+1}$
by $\bar\tau_i^{-1}$. 

The Lagrangian spheres arising from these matching
paths are well understood, and correspond to the algebraic vanishing cycles
mentioned in \fullref{rmk:lag}, spanning the orthogonal complement to
$\Lambda_2=\pi_2^*H^2(\FF_6)$ in the second homology group of $X_2$.
While it is generally expected that $X_2$ contains no ``exotic'' Lagrangian
spheres (representing homology classes that are not orthogonal to
$\Lambda_2$), the calculations in \fullref{s:mainthm} give an indication
of how one might start looking for unexpected matching paths in $\tilde{f}_2$.

For example, since the Dehn twist $\phi$ belongs to the subgroup
generated by the $\zeta_{i,j}$, $\delta_i^\pm$, and $\bar\tau_i$, one could
try to use the observation that $\zeta_{i,j}$ coincides with the conjugate
of $\xi_{i,j}$ by $(\bar\tau_1\bar\tau_2\bar\tau_3\bar\tau_4\phi)^6$ as a
starting point to build a matching path. Various other tricks of a similar
nature come to mind (all using the Dehn twist $\phi$ to build unexpected
relations among monodromy factors). However, our first naive attempts in
this direction have all led to paths that are only {\it immersed} rather
than embedded, and hence to immersed (rather than embedded) Lagrangian
spheres. The existence of such immersed spheres is  not very
surprising if one remembers that Gromov's h--principle applies to Lagrangian
immersions. In fact, the difference between immersed and embedded objects
is a recurring theme in 4--manifold topology, and it is interesting to see it
appear in this situation.

Similar considerations come up when investigating matching paths in the
canonical pencil $\tilde{f}_1$ on $X_1$. In both cases, this suggests
the following directions for further study:

\begin{question}
Can one refine the naive approach discussed above in order to
exhibit {\em embedded} ``exotic'' matching paths and Lagrangian spheres in $X_i$?

If not, what is a good way to algebraize the distinction between
``embedded'' and ``immersed'' relations among Dehn twists in a mapping
class group factorization?
\end{question}
\bigskip

\bibliographystyle{gtart}
\bibliography{link}

\begin{thebibliography}{}
\providecommand\bibmarginpar{\leavevmode\marginpar}
\def\urlstyle#1{{\tt #1}}

\bibitem{Ag2}
\textbf{D Auroux}, \emph{Fiber sums of genus 2 {L}efschetz fibrations}, Turkish
  J. Math. 27 (2003) 1--10 \xox{MR}{1975329}

\bibitem{ADK}
\textbf{D Auroux}, \textbf{S\,K Donaldson}, \textbf{L Katzarkov},
  \href{http://dx.doi.org/10.1007/s00208-003-0418-9} {\emph{Luttinger surgery
  along {L}agrangian tori and non-isotopy for singular symplectic plane
  curves}}, Math. Ann. 326 (2003) 185--203 \xox{MR}{1981618}

\bibitem{AK}
\textbf{D Auroux}, \textbf{L Katzarkov}, \emph{A degree doubling formula for
  braid monodromies and Lefschetz pencils}, preprint

\bibitem{AK1}
\textbf{D Auroux}, \textbf{L Katzarkov},
  \href{http://dx.doi.org/10.1007/PL00005795} {\emph{Branched coverings of
  $\mathbf{C}\mathrm{P}^2$ and invariants of symplectic 4--manifolds}}, Invent.
  Math. 142 (2000) 631--673 \xox{MR}{1804164}

\bibitem{AMP}
\textbf{D Auroux}, \textbf{V Mu{\~n}oz}, \textbf{F Presas}, \emph{Lagrangian
  submanifolds and {L}efschetz pencils}, J. Symplectic Geom. 3 (2005) 171--219
  \xox{MR}{2199539}

\bibitem{AS}
\textbf{D Auroux}, \textbf{I Smith}, \emph{Lefschetz pencils, branched covers
  and symplectic invariants}, from: ``Proc.\ CIME school, Symplectic
  4--manifolds and algebraic surfaces (Cetraro, 2003)'', Lecture Notes in
  Math., Springer  , to appear

\bibitem{Do1}
\textbf{S\,K Donaldson}, \emph{Lefschetz fibrations in symplectic geometry},
  from: ``Proceedings of the International Congress of Mathematicians, Vol. II
  (Berlin, 1998)'', Extra Vol. II (1998)  309--314 \xox{MR}{1648081}

\bibitem{Do2}
\textbf{S\,K Donaldson}, \emph{Lefschetz pencils on symplectic manifolds}, J.
  Differential Geom. 53 (1999) 205--236 \xox{MR}{1802722}

\bibitem{FS}
\textbf{R Fintushel}, \textbf{R\,J Stern}, \emph{Symplectic surfaces in a fixed
  homology class}, J. Differential Geom. 52 (1999) 203--222 \xox{MR}{1758295}

\bibitem{FM}
\textbf{R Friedman}, \textbf{J\,W Morgan}, \emph{Algebraic surfaces and
  {S}eiberg--{W}itten invariants}, J. Algebraic Geom. 6 (1997) 445--479
  \xox{MR}{1487223}

\bibitem{Fu}
\textbf{T Fuller}, \href{http://dx.doi.org/10.1007/s002080050182}
  {\emph{Diffeomorphism types of genus 2 {L}efschetz fibrations}}, Math. Ann.
  311 (1998) 163--176 \xox{MR}{1624287}

\bibitem{Go}
\textbf{R\,E Gompf},
  \href{http://projecteuclid.org/getRecord?id=euclid.jsg/1094072003}
  {\emph{Toward a topological characterization of symplectic manifolds}}, J.
  Symplectic Geom. 2 (2004) 177--206 \xox{MR}{2108373}

\bibitem{GS}
\textbf{R\,E Gompf}, \textbf{A\,I Stipsicz}, \emph{4--manifolds and {K}irby
  calculus}, Graduate Studies in Mathematics 20, Amer. Math. Soc. (1999)
  \xox{MR}{1707327}

\bibitem{Gr}
\textbf{M Gromov}, \emph{Pseudoholomorphic curves in symplectic manifolds},
  Invent. Math. 82 (1985) 307--347 \xox{MR}{809718}

\bibitem{HLS}
\textbf{H Hofer}, \textbf{V Lizan}, \textbf{J-C Sikorav}, \emph{On genericity
  for holomorphic curves in four-dimensional almost-complex manifolds}, J.
  Geom. Anal. 7 (1997) 149--159 \xox{MR}{1630789}

\bibitem{Hor}
\textbf{E Horikawa},
  \href{http://links.jstor.org/sici?sici=0003-486X(197609)2:104:2%3C357:ASOGT%
W%3E2.0.CO%3B2-6} {\emph{Algebraic surfaces of general type with small
  $C^{2}_{1}.$\ {I}}}, Ann. of Math. $(2)$ 104 (1976) 357--387
  \xox{MR}{0424831}

\bibitem{Kas}
\textbf{A Kas},
  \href{http://projecteuclid.org/getRecord?id=euclid.pjm/1102779371} {\emph{On
  the handlebody decomposition associated to a {L}efschetz fibration}}, Pacific
  J. Math. 89 (1980) 89--104 \xox{MR}{596919}

\bibitem{McD}
\textbf{D McDuff},
  \href{http://projecteuclid.org/getRecord?id=euclid.dmj/1092749085}
  {\emph{Almost complex structures on $S^2\times S^2$}}, Duke Math. J. 101
  (2000) 135--177 \xox{MR}{1733733}

\bibitem{McS}
\textbf{D McDuff}, \textbf{D Salamon}, \emph{$J$--holomorphic curves and
  symplectic topology}, American Mathematical Society Colloquium Publications
  52, American Mathematical Society, Providence, RI (2004) \xox{MR}{2045629}

\bibitem{MiWh}
\textbf{M\,J Micallef}, \textbf{B White},
  \href{http://links.jstor.org/sici?sici=0003-486X(199501)2:141:1%3C35:TSOBPI%%
3E2.0.CO%3B2-G} {\emph{The structure of branch points in minimal surfaces and
  in pseudoholomorphic curves}}, Ann. of Math. $(2)$ 141 (1995) 35--85
  \xox{MR}{1314031}

\bibitem{Moi}
\textbf{B\,G Moishezon}, \emph{Stable branch curves and braid monodromies},
  from: ``Algebraic geometry (Chicago, Ill., 1980)'', Lecture Notes in Math.
  862, Springer, Berlin (1981)  107--192 \xox{MR}{644819}

\bibitem{Se}
\textbf{P Seidel}, \emph{More about vanishing cycles and mutation}, from:
  ``Symplectic geometry and mirror symmetry (Seoul, 2000)'', World Sci. Publ.,
  River Edge, NJ (2001)  429--465 \xox{MR}{1882336}

\bibitem{ST1}
\textbf{B Siebert}, \textbf{G Tian},
  \href{http://dx.doi.org/10.1142/S0219199799000110} {\emph{On hyperelliptic
  $C^\infty$--{L}efschetz fibrations of four-manifolds}}, Commun. Contemp.
  Math. 1 (1999) 255--280 \xox{MR}{1696101}

\bibitem{ST}
\textbf{B Siebert}, \textbf{G Tian}, \emph{On the holomorphicity of genus two
  {L}efschetz fibrations}, Ann. of Math. $(2)$ 161 (2005) 959--1020
  \xox{MR}{2153404}

\bibitem{Sik}
\textbf{J-C Sikorav}, \emph{The gluing construction for normally generic
  $J$--holomorphic curves}, from: ``Symplectic and contact topology:
  interactions and perspectives (Toronto, ON/Montreal, QC, 2001)'', Fields
  Inst. Commun. 35, Amer. Math. Soc., Providence, RI (2003)  175--199
  \xox{MR}{1969276}

\bibitem{Sm}
\textbf{I Smith}, \href{http://dx.doi.org/10.2140/gt.2001.5.579}
  {\emph{Lefschetz pencils and divisors in moduli space}}, Geom. Topol. 5
  (2001) 579--608 \xox{MR}{1833754}

\end{thebibliography}

\end{document}